\documentclass[11pt]{article}

\usepackage{natbib}
\usepackage{amsmath,amssymb}
\usepackage{color} 
\usepackage{graphicx}
\def\ip<#1,#2>{\left\langle #1,#2\right\rangle}

\newcommand\Rey{\mbox{\textit{Re}}}  % Reynolds number
\newsavebox{\astrutbox}
\sbox{\astrutbox}{\rule[-5pt]{0pt}{20pt}}

\def\drawline#1#2{\raise 2.5pt\vbox{\hrule width #1pt height #2pt}}
\def\spacce#1{\hskip #1pt}
\def\solid{\drawline{24}{.5}\nobreak}
\def\bdash{\hbox{\drawline{4}{.5}\spacce{2}}}
\def\dashed{\bdash\bdash\bdash\bdash\hskip-2pt\nobreak}
\def\bdot{\hbox{\drawline{1}{.5}\spacce{2}}}
\def\dotted{\hbox{\leaders\bdot\hskip 24pt}\hskip-2pt\nobreak}
\def\chndash{\hbox {\drawline{8.5}{.5}\spacce{2}\drawline{3}{.5}\spacce{2}\drawline{8.5}{.5}}\nobreak}
\def\chndot{\hbox {\drawline{9.5}{.5}\spacce{2}\drawline{1}{.5}\spacce{2}\drawline{9.5}{.5}}\nobreak}

\def\trian{\raise 1.25pt\hbox{$\scriptscriptstyle\triangle$}\nobreak}
\def\solidtrian{\raise 1.25pt
\hbox to 3bp{
\hfill}\nobreak}

\def\plus{\raise 1.25pt \hbox{$\scriptscriptstyle +$}\nobreak\ }
\def\x{\raise 1.25pt \hbox{$\scriptscriptstyle \times$}\nobreak\ }

\graphicspath{{arxiv_figs/}}
\begin{document}

\title{Feedback control of unstable steady states of flow past a flat plate using reduced-order estimators}

%\author[S. Ahuja and C. W. Rowley]%
%{S.\ns A\ls H\ls U\ls J\ls A $^1$ \ns%
%\and C.\ns W.\ns R\ls O\ls W\ls L\ls E\ls Y $^1$}

%\affiliation{$^1$Department of Mechanical and  Aerospace Engineering,\\ Princeton University, Princeton, NJ 08544, USA}

\author{Sunil Ahuja and Clarence W. Rowley \\
\small{\{sahuja, cwrowley\} at princeton.edu} \\
\small{\em Mechanical and Aerospace Engineering,} \\
\small {\em Princeton University, Princeton, NJ 08544, USA.}
}

\maketitle \thispagestyle{empty}

\begin{abstract}

We present an estimator-based control design procedure for flow control, using reduced-order models of the governing equations, linearized about a possibly unstable steady state. The reduced-order models are obtained using an approximate balanced truncation method that retains the most controllable and observable modes of the system. The original method is valid only for stable linear systems, and in this paper, we present an extension to unstable linear systems. The dynamics on the unstable subspace are represented by projecting the original equations onto the global unstable eigenmodes, assumed to be small in number. A snapshot-based algorithm is developed, using approximate balanced truncation, for obtaining a reduced-order model of the dynamics on the stable subspace.

The proposed algorithm is used to study feedback control of two-dimensional flow over a flat plate at a low Reynolds number and at large angles of attack, where the natural flow is vortex shedding, though there also exists an unstable steady state. For control design, we derive reduced-order models valid in the neighborhood of this unstable steady state. The actuation is modeled as a localized body force near the leading edge of the flat plate, and the sensors are two velocity measurements in the near-wake of the plate. A reduced-order Kalman filter is developed based on these models and is shown to accurately reconstruct the flow field from the sensor measurements, and the resulting estimator-based control is shown to stabilize the unstable steady state. For small perturbations of the steady state, the model accurately predicts the response of the full simulation. Furthermore, the resulting controller is even able to suppress the stable periodic vortex shedding, where the nonlinear effects are strong, thus implying a large domain of attraction of the stabilized steady state.

\end{abstract}

\section{Introduction}\label{sec:introduction}

The goal of this paper is two-fold; the first goal is to present an algorithm for developing reduced-order models of the input-output dynamics of {\em unstable} high-dimensional linear state-space systems (such as linearized Navier-Stokes equations with actuation and sensing), while the second goal is to demonstrate the algorithm by developing estimation-based controllers to stabilize unstable steady states of a two dimensional low-Reynolds-number flow past a flat plate at a large angle of attack.

\subsection{Model reduction for unstable systems}

Development of feedback control strategies based on linearized Navier-Stokes equations is attractive due to the ready availability of a large class of control techniques, and there has been substantial progress in this direction in the past decade, reviewed in detail by~\cite{KimBew-07}. However, many of these techniques are limited to relatively small dimensional systems~$\sim O(10^3)$, while the numerical discretization of fluid flows invariably result in huge dimensional systems, typically~$O(10^{5-8})$. Thus, model reduction has played an important role in making these tools further accessible to fluid flows.

Extensive research effort in model reduction has focused on the method of proper orthogonal decomposition~(POD) and Galerkin projection, developed first by~\cite{Lumley}. The main disadvantage of this technique is that, although the POD modes capture the energetically important structures of the flow, the reduced-order models obtained by the subsequent Galerkin projection of the governing equations onto these modes often do not faithfully represent the dynamics. Various modifications to improve this method have been proposed and used for flow control; refer to the introduction of~\cite{Siegel-jfm08} for a review of these techniques. The POD/Galerkin methods have been applied for flow control in various contexts, such as bluff-body wake suppression (\cite{Siegel-jfm08}, \cite{GrPeTa-99}, \cite{Noack:2004}, \cite{TadmorCNLLM-07}), noise reduction in cavity flow~(\cite{Rowley-cdc05}, \cite{Gloerfelt-08}), and drag reduction in turbulent boundary layers~(\cite{LumleyBlossey-98}, \cite{PrabhuCoCh-01}). Another model-reduction technique, based on projection onto the global (stable or unstable) eigenmodes of the flow linearized about steady states, has been used by~\cite{Aakervik_et_al-jfm-07} and~\cite{HenAk08} in the context of spatially developing flows such as separated boundary layers.
% A disadvantage of this procedure is that it does not take the actuation or sensing into account, which can be important for flow control.
% The technique works very well for spatially-developing wall-flows as the linearized and adjoint modes have support near the wall, feasible for actuator and sensor placement. For wakes, we will see that the unstable eigenmodes have little or no support close to the wall, making the sensor problem much more difficult. (Perhaps mention this point later in the section on eigenmodes?)
% In flows dominated by vortices, vortex-models have been developed and used for control by {Cortellezi96}
In this paper, we focus on an approximate balanced truncation method developed by~\cite{Rowley-ijbc05} as an approximation to the original method of~\cite{Moore-81}. This technique captures the dynamically important modes of the system, and the non-approximate version provides rigorous bounds for the resulting reduced-order models. The method, sometimes called balanced POD, was used to obtain models of the linearized channel flow by~\cite{IlakRowley-pof08} and the Blasius boundary layer by~\cite{BaBrHe-09}, and was shown to accurately capture the control actuation and also to outperform the POD/Galerkin models.

The balanced truncation method of~\cite{Moore-81} is applicable only to systems linearized about {\em stable} steady states. An extension to {\em unstable} linear systems was proposed by~\cite{ZhoSalWu-99}, by introducing frequency-domain definitions of controllability and observability Gramians. Reduced-order models were obtained by first decoupling the dynamics on the stable and unstable subspaces, and then truncating the relatively uncontrollable and unobservable modes on each of the two subspaces. In this paper, we present an approximation algorithm for balanced truncation of linear unstable systems, which results in models that are equivalent to those of~\cite{ZhoSalWu-99} on the {\em stable} subspace. The dynamics on the unstable subspace is treated {\em exactly} by a projection onto the global eigenmodes, as in~\cite{Aakervik_et_al-jfm-07}.

% Models using unstable global eigenmodes for spatially developing flow is feasible, as the linearized and adjoint modes both have support near the wall, implying feasible sensor and actuator locations. However, for wake flows, the problem of sensing is much harder, as the linearized unstable modes have support at a distance downstream from the wall, while the ideal sensors would perhaps be pressure measurements on the surface of the plate.

\subsection{Control of flow past 2-D wings}

As a proof-of-concept study, the modeling procedure is applied to the problem of two dimensional low-Reynolds-number flow past a flat plate at a large angle of attack. We develop reduced-order models and design controllers that stabilize the unstable steady states of this flow. Our motivation for the choice of this problem comes from our interest in regulating vortices in separated flows behind low aspect-ratio wings, which is of importance in design of micro air vehicles~(MAVs). Recently, design of MAVs has been inspired from experimental observations in insect and bird flights of a stabilizing leading edge vortex (see~\cite{BirDic-01} and~\cite{EllingtonVaWiTh-96}), which remains attached throughout the wing stroke and provides enhanced lift. So, it could be beneficial to design controllers that can manipulate the wake of MAVs to enhance lift and achieve better maneuverability in presence of wind gusts. Recent studies in this direction, using open-loop control of the flow past low-aspect-ratio wings using steady or periodic blowing, were performed computationally  by~\cite{TairaCol-aiaaj09} and experimentally by~\cite{WillCJCT-08}. These studies explored different forcing amplitudes and frequencies, locations and directions. However, the design of feedback controllers remains a challenge, due to the large dimensionality of the problem and complex flow physics. We present computational tools, that we hope can at least pave a direction and provide techniques towards addressing some of these challenges. We consider the 2-D flow past a flat plate, actuated by a localized body force close to the leading edge, with two near-wake velocity sensors. We design a reduced-order compensator and show that it is able to suppress vortex shedding at high angles of attack.

%In experiments on rapidly pitching airfoils, much higher lift is observed while the angle of attack is increasing as compared to decreasing, and is attributed to the formation of a dynamic stall vortex. This vortex sheds when the airfoil pitches down, resulting in a drop in lift. Recent experiments by~\cite{WillCJCT-08} have been able to achieve the high lift branch observed during pitch-up even for slowly pitching using open-loop control. For stationary airfoils translating at low speeds, as the angle of attack is increased, the flow undergoes a Hopf bifurcation from a steady state to periodic vortex shedding. However, an {\em unstable} steady state also exists at these high angles of attack, and could be {\em stabilized} using feedback control. For 2-D flows past a flat plate, we find that the lift corresponding to these unstable steady states is close to the minimum of that in the periodic vortex shedding cycle (see~Fig.~\ref{fig:ss}). However, the flow physics of low aspect ratio flat plate is considerably different, due to strong interaction of the near-wake with tip vortices, as reported by~\cite{TairaCol-jfm08}. In this paper, we focus on computational tools to develop controllers that stabilize unstable steady states in a 2-D flow past a flat plate. Future studies will investigate the existence of high-lift unstable steady states in flow past low aspect ratio wings, and focus on development of controllers to achieve high-lift flow states in steady MAV flight. {\color{blue} (TODO: Work more here.) }

Many previous studies have focused on the control of flow past a cylinder, which is qualitatively similar to the flow past a flat plate at large angle of attack, with the natural flow in both the cases being periodic vortex shedding. For the cylinder, the flow undergoes a transition from steady state to periodic shedding with increasing Reynolds number, while a similar transition occurs in the flat plate with increasing angle of attack. There has been considerable research effort on suppression of this shedding in cylinder and other bluff body wakes, using passive and active, open-loop and feedback control, as reviewed by~\cite{ChoiJeonKim08}. Among those, some techniques are based on reduced-order models; for instance, \cite{Gillies-98} developed models using artificial neural networks and a POD basis, \cite{GrPeTa-99} modified the POD/Galerkin method to account for actuation by means of cylinder-rotation, while~\cite{Siegel-jfm08} developed a double POD method to account for changes in the wake structure during transients. Some earlier efforts in the control of a {\em flat-plate} wake include those by~\cite{Cortelezzi-96}, \cite{CorCheCha-97} who used vortex-based methods to model the flow past a vertical plate~(angle of attack~=~$90^\circ$); vortex-based models form their own class of modeling techniques reviewed recently by~\cite{Protas08}. Lagrangian coherent structures were used by~\cite{WaHaBaTadmor-03} to enhance mixing in flow past a bluff body with the trailing surface similar to the vertical flat plate. One of the few efforts towards control of flat plate at an angle of incidence was by \cite{ZaIollo-03}, who used a passive leading-edge suction control along with a potential flow vortex model. \cite{PaHeNoKiTad-jfm08} also used reduced-order vortex models for drag reduction on an elongated D-shaped bluff body.

% Control of flow past a plate at right angles to the flow direction was studied by \cite{Iollo03,Cortellezi96,Cortellezi97,Shermer92}, \cite{Sarpakaya75, KiyaArie77}
% Our aim here is to develop techniques that will eventually be used towards our goal of being able to manipulate vortices over models real wings, and be integrated with flight control
% The tools that we present here are a lot more general, and can be extended to other flow control scenarios

% The standard benchmark for control strategies has been that of flow past a cylinder, however, due to our interest in gradually moving towards real wings, we consider the flow past a flat plate at large angles of attack
% Most of the work on control of flow past airfoils has been experimental, and there are relatively fewer reports on model-based control strategies.

This paper is organized as follows: In section~\ref{sec:model_red}, we first briefly describe the balanced truncation method for unstable systems as developed by~\cite{ZhoSalWu-99}, and the approximate balanced truncation procedure called balanced POD of~\cite{Rowley-ijbc05} for large dimensional stable systems. Then, we present an algorithm for approximate balanced truncation of large dimensional {\em unstable} systems, assuming that the dimension of unstable subspace is small and the corresponding global eigenmodes can be computed. In section~\ref{sec:numerics_ibfs}, we briefly describe the numerical technique of~\cite{ColTai-07} using a fast immersed boundary method, and present the linearized and adjoint formulations of this numerical method. In section~\ref{section:model_results}, we present numerical results, using the model example of two-dimensional flow past a flat plate at a large angle of attack and a low Reynolds number~$\Rey = 100$. First, we perform a steady-state analysis and compute the branch of steady states in the entire range of angles of attack,~$0 \leq \alpha \leq 90^\circ$. We also compute the left and right global eigenmodes of the flow linearized about the unstable steady state at~$\alpha = 35^\circ$. We present reduced-order models of the linearized dynamics and use linear optimal control techniques to design controllers that stabilize this steady state. Full-state feedback and more practical near-wake velocity-measurement based feedback controllers are derived, implemented in the nonlinear equations, and shown to suppress vortex shedding. The paper concludes with a brief discussion in section~\ref{sec:discussion}.

\section{Model reduction methodology} \label{sec:model_red}

\subsection{Balanced truncation of unstable systems}
\label{sec:unstable_baltrunc}

We briefly describe a model reduction procedure using the balanced truncation method for unstable systems developed by~\cite{ZhoSalWu-99}. Consider the state-space system
\begin{align}
\dot{x} &= Ax + Bu, \nonumber \\
y &= Cx, \label{ss}
\end{align}
where~$x\in \mathbb{R}^n$ is the state, $u \in \mathbb{R}^p$~is the input, and $y \in \mathbb{R}^q$ is the output of the system; the dot over~$x$ represents differentiation with respect to time. The eigenvalues of~$A$ are assumed to be anywhere on the complex plane, except on the imaginary axis.

The standard balanced truncation procedure developed by~\cite{Moore-81}, valid only for stable systems, starts with defining the controllability and observability Gramians of the system~(\ref{ss}) as follows:
\begin{align}
W_c &= \int_0^\infty e^{At} B B^\ast e^{A^\ast t} \, dt \nonumber \\
\mbox{and} \quad W_o &= \int_0^\infty e^{A^\ast t} C^\ast C e^{A t} \, dt,
\label{gramians_stable}
\end{align}
where asterisks are used to denote adjoint operators. A co-ordinate transformation is then obtained such that the Gramians~(\ref{gramians_stable}) of the transformed system are equal and diagonal. The diagonal entries of the transformed Gramians, called Hankel singular values~(HSVs), decrease monotonically and are directly related to the controllability and observability of the corresponding states. A reduced-order model is obtained by truncating the states with relatively small HSVs, that is, the states which are almost uncontrollable and unobservable.

For unstable systems, the integrals in~(\ref{gramians_stable}) are unbounded and hence the Gramians are ill-defined. A modified technique was proposed by~\cite{ZhoSalWu-99} based on the following frequency-domain definitions of the Gramians:
\begin{align}
W_c &= \frac{1}{2\pi} \int_{-\infty}^{\infty} (j\omega I -A)^{-1} B B^\ast (-j\omega I - A^\ast)^{-1} \, \, d\omega,  \\
W_o &= \frac{1}{2\pi} \int_{-\infty}^{\infty} (-j\omega I -A^\ast)^{-1} C^\ast C (j\omega I - A)^{-1} \, \, d\omega.
\label{gram_freq}
\end{align}
By using Parseval's theorem, it can be shown that for stable systems, the frequency-domain definitions~(\ref{gram_freq}) are equivalent to the time-domain definitions~(\ref{gramians_stable}). The model-reduction procedure of~\cite{ZhoSalWu-99} begins by first transforming the system~(\ref{ss}) to coordinates in which the stable and unstable dynamics are decoupled. That is, let $T$~be a transformation such that if~$x = T \tilde{x}$, the
system~(\ref{ss}) transforms to
\begin{align}
\dot{\tilde{x}} = \frac{d}{dt} \begin{pmatrix} \tilde{x}_u \\ \tilde{x}_s \end{pmatrix} & =
\begin{pmatrix} A_u & 0 \\ 0 & A_s \end{pmatrix} \tilde{x} +
\begin{pmatrix} B_u \\ B_s \\ \end{pmatrix} u \nonumber \\
y & = \begin{pmatrix} C_u & C_s  \end{pmatrix} \tilde{x}.
\label{decouple_ss}
\end{align}
Here, $A_u$~and~$A_s$ are such that all their eigenvalues are in the right- and left-half complex planes respectively, while~$\tilde{x}_u$ and~$\tilde{x}_s$ are the corresponding states. Next, denote the controllability and observability Gramians corresponding to the set~$(A_s, B_s, C_s)$ describing the stable dynamics by $W_c^s$~and~$W_o^s$ respectively. Similarly, denote the Gramians corresponding to the set~$(-A_u, B_u, C_u)$ by~$W_c^u$~and~$W_o^u$. The Gramians of the original system are then related to those corresponding to the subsystems by:
\begin{align}
W_c & = T \begin{pmatrix} W_c^u & 0 \\ 0 & W_c^s \end{pmatrix} T^\ast \nonumber \\
\mbox{and} \quad W_o & = (T^{-1})^\ast \begin{pmatrix} W_u^u & 0 \\ 0 & W_u^s \end{pmatrix} T^{-1}. \label{gramians_general}
\end{align}
A system is said to be balanced if its Gramians defined by~(\ref{gramians_general}) are equal and diagonal, in which case the diagonal entries are called the {\em generalized} Hankel singular values. A reduced-order model is obtained by truncating the states with small generalized HSVs.

\subsection{Approximate balanced truncation of stable systems}
\label{sec:approx_baltrunc_stable}

For systems of large dimension~$\sim O(10^{5-8})$ such as those encountered here, the Gramians~(\ref{gramians_general}) are huge matrices which cannot be easily computed or stored. A computationally tractable procedure was introduced by~\cite{Rowley-ijbc05} for obtaining an approximate balancing transformation. We first briefly describe this method, valid only for stable systems, and then present an extension to unstable systems. The procedure consists of computing the impulse responses of the system~(\ref{ss}) and stacking the resulting snapshots of the state~$x$ as columns of a matrix~$X$. It also requires state-snapshots of the impulse responses of the adjoint system
\begin{align}
\dot{z} & = A^\ast z + C^\ast v, \nonumber \\
w &= B^\ast z, \label{ss_adj}
\end{align}
which are stacked as columns of a matrix~$Z$. Then, the Gramians~(\ref{gramians_stable}) can be approximated as
\begin{align}
W_c \approx X X^\ast, \qquad W_o \approx Z Z^\ast.
\label{gramians_approx}
\end{align}
The approximate Gramians~(\ref{gramians_approx}) are not actually computed due to the large storage cost, but the leading columns (or modes) of the transformation that balances these Gramians are computed using a cost-efficient algorithm. It involves computing the singular value decomposition of
\begin{align}
Z^\ast X = U \Sigma V^\ast =
\begin{pmatrix} U_1 & U_2 \\ \end{pmatrix}
\begin{pmatrix} \Sigma_1 & 0 \\ 0 & \Sigma_2 \end{pmatrix}
\begin{pmatrix} V_1^\ast \\ V_2^\ast \\ \end{pmatrix},
\label{svd}
\end{align}
where~$\Sigma_1 \in \mathbb{R}^{r\times r}$ is a diagonal matrix of the most significant HSVs greater than a cut-off which is a modeling parameter, while~$\Sigma_2 \in \mathbb{R}^{(n-r)\times (n-r)}$ is a diagonal matrix of smaller and zero HSVs. Note that $Z^\ast X \in \mathbb{R}^{n_o \times n_s}$ is a small matrix, where~$n_s$ and~$n_o$ are the number of snapshots of the impulse responses of systems~(\ref{ss}) and~(\ref{ss_adj}) respectively. For fluid systems that we are interested in, the typical number of snapshots is~$O(10^{2-4})$, thus resulting in a reasonable computational cost, and typically~$r \leq 50$. The leading columns and rows of the balancing transformation and its inverse are obtained using:
\begin{align}
\Phi = XV_1\Sigma_1^{-1/2}, \quad  \Psi = ZU_1\Sigma_1^{-1/2},
\label{bal_transf}
\end{align}
where~$\Phi, \Psi \in \mathbb{R}^{n \times r}$, and the two sets of modes are bi-orthogonal; that is,~$\Psi^\ast\Phi=I$. The reduced-order model of~(\ref{ss}) is then obtained by expressing~$x = \Phi a$, $a \in \mathbb{R}^r$, and using the bi-orthogonality of~$\Phi$~and~$\Psi$:
\begin{align}
\dot{a} & = \Psi^\ast A \Phi a + \Psi^\ast B u, \\
y & = C \Phi a.
\end{align}

\subsubsection{Output projection}
\label{sec:outproj}

When the number of outputs of the system (rows of~$C$) is large, the algorithm described in section~\ref{sec:approx_baltrunc_stable} can become intractable. The reason for this is that it involves one simulation of the adjoint system~(\ref{ss_adj}) for each component of~$v$, the dimension of which is the same as the number of outputs. This number is often large in fluids systems where a good model needs to capture the response of the entire system~($C=I$) to a given input. To overcome this problem,~\cite{Rowley-ijbc05} proposed a technique called~{\em output projection}, which involves projecting the output~$y$ of~(\ref{ss}) onto a small number of energetically important modes obtained using proper orthogonal decomposition~(POD). Let the columns of~$\Theta \in \mathbb{R}^{q \times m}$ consist of the first $m$~POD modes of the dataset consisting of outputs obtained from an impulse response of~(\ref{ss}). Then, for the purpose of obtaining a reduced-order model, the output of~(\ref{ss}) is approximated by
\begin{align}
y &= \Theta \Theta^\ast C x, \label{ss_outproj}
\end{align}
where~$\Theta \Theta^\ast$ is an orthogonal projection of the output onto the first $m$~POD modes. The resulting output-projected system is optimally close (in the $L^2-$sense) to the original system, for an output of fixed rank~$m$. With this approximation, only~$m$ adjoint simulations are required to approximate the observability Gramian; refer to~\cite{Rowley-ijbc05} for details. The number of POD modes~$m$ for output projection is a design parameter and is typically chosen to capture more than $90\%$ of the output energy. In the rest of this paper, the models resulting from this approximation of the output are referred to as {\em $m-$mode output projected models.}

%The reduced-order model of the output projected system is then given by
%%
%\begin{align}
%\dot{a} & = \Psi^\ast A \Phi a + \Psi^\ast B u, \\
%y & = \Theta \Theta^\ast C \Phi a.
%\end{align}

\subsection{Approximate balanced truncation of unstable systems}
\label{sec:approx_baltrunc_unstable}

The approximate balancing procedure described in the previous section, which is essentially a snapshot-based method, does not extend to unstable systems since the impulse responses of~(\ref{ss}) and~(\ref{ss_adj}) are unbounded. We could consider applying the algorithm to the two sub-systems given in~(\ref{decouple_ss}), but the transformation~$T$ that decouples~(\ref{ss}) itself is not available. However, when the dimension of the unstable sub-system is small, we show that it is not necessary to compute the entire transformation~$T$ and it is still possible to obtain an approximate balancing transformation. Here, we present an algorithm for computing such a transformation and also show that it essentially results in a method that is a slight variant of the technique of~\cite{ZhoSalWu-99} presented in section~\ref{sec:unstable_baltrunc}. The idea behind the algorithm is to project the original system~(\ref{ss}) onto the stable subspace of~$A$. Then, one obtains a reduced-order model of the projected system using the snapshot-based procedure described in section~\ref{sec:approx_baltrunc_stable}. The dynamics projected onto the unstable subspace can be treated exactly on account of its low dimensionality.

We assume that the number of unstable eigenvalues~$n_u$ is~$O(10)$ and can be computed numerically, say using the computational package ARPACK developed by~\cite{arpack98}. We further assume that the bases for the right and the left unstable eigenspaces~$\Phi_u, \Psi_u \in \mathbb{R}^{n \times n_u}$ can be computed. For the algorithm, we need the following projection operator onto the stable subspace of~$A$:
\begin{align}
\mathbb{P}_s = I - \Phi_u \Psi_u^\ast, \label{proj_stable}
\end{align}
where~$\Phi_u$ and~$\Psi_u$ have been scaled such that~$\Psi_u^\ast \Phi_u=I_{n_u}$. We use the operator~$\mathbb{P}_s$ to obtain the dynamics of~(\ref{ss}) restricted to the stable subspace of~$A$ as follows:
\begin{align}
\dot{x}_s &=  A x_s + \mathbb{P}_s B u, \nonumber \\
y_s &= C \mathbb{P}_s x_s  \label{ss_stable}
\end{align}
where~$x_s = \mathbb{P}_s x.$ The adjoint of~(\ref{ss_stable}) is the same as the dynamics of~(\ref{ss_adj}) restricted to the stable subspace of~$A^\ast$ using~$\mathbb{P}_s^\ast$, and is given by
\begin{align}
\dot{z}_s & = A^\ast z_s + \mathbb{P}_s^\ast C^\ast v,  \nonumber \\
w_s & = B^\ast \mathbb{P}_s^\ast z_s, \label{ss_adj_stable}
\end{align}
where~$z_s = \mathbb{P}_s^\ast z$. We compute the state-impulse responses of~(\ref{ss_stable}) and~(\ref{ss_adj_stable}) and stack the resulting snapshots~$x_s$ and~$z_s$ in matrices~$X_s$ and~$Z_s$ respectively. As in~(\ref{svd}), we compute the singular valued decomposition of~$Z_s^\ast X_s$ and use the expressions~(\ref{bal_transf}) to obtain the balancing modes~$\Phi_s$ and the adjoint modes~$\Psi_s$, where again~$\Psi_s^\ast \Phi_s = I_r$. The reduced-order modes are obtained by expressing the state~$x$ as
\begin{align}
x = \Phi_u a_u + \Phi_s a_s, \label{modal_exp}
\end{align}
where~$a_u \in \mathbb{R}^{n_u}$ and~$a_s \in \mathbb{R}^{r}$. Substituting~(\ref{modal_exp}) in~(\ref{ss}) and pre-multiplying by~$\Psi_u^\ast$ and~$\Psi_s^\ast$, we obtain
\begin{align}
\frac{da}{dt} \equiv \frac{d}{dt} \begin{pmatrix} a_u \\ a_s \\ \end{pmatrix} & =
\begin{pmatrix}  \Psi_u^\ast A \Phi_u & \Psi_u^\ast A \Phi_s \\ \Psi_s^\ast A \Phi_u & \Psi_s^\ast A \Phi_s \\ \end{pmatrix}
\begin{pmatrix} a_u \\ a_s \\ \end{pmatrix} +
\begin{pmatrix}  \Psi_u^\ast \\ \Psi_s^\ast \\ \end{pmatrix} Bu \label{romA}\\
y & = C (\Phi_u a_u + \Phi_s a_s) \equiv \begin{pmatrix} C\Phi_u & C\Phi_s \end{pmatrix} a. \label{romB}
\end{align}
Now, since~range$(A \Phi_u) \subseteq$~span$(\Phi_u)$, we can write~$A \Phi_u = \Phi_u \Lambda$ for some~$\Lambda \in \mathbb{R}^{n_u \times n_u}$, and using the properties of eigenvectors, we have~$\Psi_s^\ast A \Phi_u = \Psi_s^\ast \Phi_u \Lambda = 0$. Similarly, it can be shown that~$\Psi_u^\ast A \Phi_s = 0$. Thus, the cross terms in~(\ref{romA}) are zero and the reduced-order model is
\begin{align}
\frac{da}{dt} & =
\begin{pmatrix}  \Psi_u^\ast A \Phi_u & 0 \\ 0 & \Psi_s^\ast A \Phi_s \\ \end{pmatrix}
\begin{pmatrix} a_u \\ a_s \\ \end{pmatrix} +
\begin{pmatrix}  \Psi_u^\ast \\ \Psi_s^\ast \\ \end{pmatrix} Bu
\equiv \begin{pmatrix} \widetilde{A}_u & 0 \\ 0 & \widetilde{A}_s \end{pmatrix}
\begin{pmatrix} a_u \\ a_s \\ \end{pmatrix} +
\begin{pmatrix} \widetilde{B}_u \\ \widetilde{B}_s \end{pmatrix} u \nonumber\\
y & = C (\Phi_u a_u + \Phi_s a_s) \equiv \begin{pmatrix} \widetilde{C}_u & \widetilde{C}_s\end{pmatrix} a.
\label{rom}
\end{align}
The procedure described so far to obtain the reduced-order model~(\ref{rom}) is related to the procedure of~\cite{ZhoSalWu-99} described in section~\ref{sec:unstable_baltrunc}. It can be shown that the transformation that balances the Gramians defined by~(\ref{gramians_general}) results in a system in which the unstable and stable dynamics are decoupled. Furthermore, the resulting stable dynamics are the same as those given by the equations describing the $a_s$-dynamics of~(\ref{rom}). That is, balancing the stable part of the Gramians~$W_c$ and~$W_o$ defined in~(\ref{gramians_general}) (balancing~$W_c^s$ and~$W_o^s$)~is the same as balancing the Gramians of the stable subsystem~(\ref{ss_stable}); a proof is outlined in appendix~\ref{sec:app:baltrunc}. In our algorithm, the unstable dynamics are not balanced, while they are by~\cite{ZhoSalWu-99}. A disadvantage of Zhou's approach is that an unstable mode might be truncated resulting in a model which does not capture the instability, which is undesirable for control purposes.

\subsubsection{Output projection for the stable subspace} \label{sec:outproj_unstable}

For systems with a large number of outputs, the number of adjoint simulations~(\ref{ss_adj_stable}) can become intractable; however, the output projection of section~\ref{sec:outproj} can readily be extended to unstable systems. Instead of projecting the entire output~$y$ onto POD modes, we first express the state~$x = x_u + x_s$, where~$x_u = (I-\mathbb{P}_s)x$ and~$x_s = \mathbb{P}_s x$ are projections on the unstable and stable subspaces of~$A$ respectively. We similarly express the output as~$y = y_u + y_s = C(x_u + x_s)$. We then project the component~$y_s$ onto a small number of POD modes, of the data set consisting of the outputs from an impulse response of~(\ref{ss_stable}). If the POD modes are represented as columns of the matrix~$\Theta_s \in \mathbb{R}^{q \times m}$, the output of~(\ref{ss}) is approximated by
\begin{align}
y = \big[ C(I - \mathbb{P}_s) + \Theta_s \Theta_s^\ast C \mathbb{P}_s \big] x = Cx_u + \Theta_s \Theta_s^\ast C x_s.
\label{yapprox_unstable}
\end{align}
Finally, with the state~$x$ expressed by the modal expansion~(\ref{modal_exp}), the output of the reduced-order model~(\ref{rom}) is given by
\begin{align}
y = \begin{pmatrix} C\Phi_u &  \Theta_s \Theta_s^\ast C \Phi_s \end{pmatrix} \begin{pmatrix} a_u \\ a_s \end{pmatrix}.
\label{yrom_unstable_outproj}
\end{align}

\subsubsection{Algorithm} \label{sec:algorithm}
The steps involved in obtaining reduced-order models of~(\ref{ss}), for the case where the output is the entire state~($C=I$), can now be summarized as follows:
\begin{enumerate}
\item \label{primal_step}
Project the original system~(\ref{ss}) onto the subspace spanned by the stable eigenvectors of~$A$ in the direction of the unstable eigenvectors of~$A$ to obtain~(\ref{ss_stable}). Compute the (state) response to an impulse on each input of~(\ref{ss_stable}) and stack the snapshots in a matrix~$X_s$.
\item
Assemble the resulting snapshots, and compute the POD modes~$\theta_j$ of the resulting data-set. These POD modes are stacked as columns of~$\Theta_s$.
\item
Choose the number of POD modes one wants to use to describe the output of~(\ref{ss_stable}).  For instance, if $10\%$ error is acceptable, and the first~$m$~POD modes capture $90\%$ of the energy, then the output is the velocity field projected onto the first $m$~modes. Thus, the output is represented as~$y_s = \Theta_s^\ast x_s$.
\item \label{adjoint_step}
Project the adjoint system~(\ref{ss_adj}) onto the subspace spanned by the stable eigenvectors of~$A^\ast$ in the direction of the unstable eigenvectors of~$A^\ast$ to obtain~(\ref{ss_adj_stable}). Compute the (state) response of~(\ref{ss_adj_stable}), starting with each POD mode~$\theta_j$ as the initial condition (one simulation for each of the first $m$~modes). Stack the snapshots in a matrix~$Z_s$.
\item
Compute the singular value decomposition of~$M = Z_s^\ast X^\ast =U_s \Sigma_s V_s^\ast$, where~$\Sigma_s \in \mathbb{R}^{r\times r}$, and~$r =$ rank($M$).
\item
Define balancing modes~$\varphi^s_j$ and the corresponding adjoint modes~$\psi^s_j$ as columns of the matrices~$\Phi_s$ and~$\Psi_s$, where
\begin{equation}
\Phi_s = X_s V_s \Sigma_s^{-1/2},\qquad \Psi_s = Y_s U_s \Sigma_s^{-1/2}. \label{phipsi_def}
\end{equation}
\item
Obtain the reduced-order model using~(\ref{rom}), which can be written as
\begin{align}
\frac{da}{dt} & = \begin{pmatrix} \widetilde{A}_u & 0 \\ 0 & \widetilde{A}_s \end{pmatrix} a +
\begin{pmatrix} \widetilde{B}_u \\ \widetilde{B}_s \end{pmatrix} u
\equiv \widetilde{A} a + \widetilde{B} u, \label{rom_compact} \\
y &= \begin{pmatrix} \widetilde{C}_u &  \widetilde{C}_s \end{pmatrix}  a \equiv \widetilde{C}a \qquad \mbox{where,} \label{output_compact}\\
a & = \begin{pmatrix} a_u \\ a_s \\ \end{pmatrix}, \\
\widetilde{A}_u & = \Psi_u^\ast A \Phi_u, \quad \widetilde{B}_u =  \Psi_u^\ast B, \quad \widetilde{C}_u =   \Phi_u, \\
\widetilde{A}_s & = \Psi_s^\ast A \Phi_s, \quad \widetilde{B}_s = \Psi_s^\ast B, \quad \widetilde{C}_s = \Theta_s \Theta_s^\ast \Phi_s.
\end{align}
Alternatively, the outputs can be considered to be simply the coefficients of the unstable modes~$a_u$ and the coefficients of the POD modes~$\Theta_s$ of the stable subspace. With this choice, the output can be represented as
\begin{align}
y &= \begin{pmatrix} \widehat{C}_u & 0 \\ 0 &  \widehat{C}_s \end{pmatrix}  \begin{pmatrix} a_u \\ a_s \end{pmatrix} \equiv \widehat{C} a, \qquad \mbox{where,}
\label{output_compact2} \\
\widehat{C}_u &= I_{n_u}, \quad \widehat{C}_s = \Theta_s^\ast \Phi_s.
\end{align}
Finally, if the initial state~$x_0$ is known, the initial condition of~(\ref{rom_compact}) can be obtained using
\begin{align}
a_0 = \begin{pmatrix} \Psi_u & \Psi_s\end{pmatrix}^\ast x_0.
\label{init_model}
\end{align}
\end{enumerate}

In the remainder of this paper, we demonstrate the algorithm developed in this section by obtaining reduced order models of the 2-D uniform flow past a flat plate, and develop controllers based on these models to stabilize the unstable steady states that exist at high angles of attack.

\section{Numerical scheme}\label{sec:numerics_ibfs} % (fold)

The numerical scheme used is a fast immersed boundary method developed by~\cite{ColTai-07}, and is briefly described here. The method is then used to develop the linearized and adjoint formulations. Consider the following form of the incompressible Navier-Stokes equations, based on the continuous analog of the immersed boundary formulation introduced by~\cite{Peskin-72}:
\begin{align}
\partial_t u + u \cdot \nabla u &= - \nabla p + \frac{1}{\Rey} \nabla^2 u + \int f(\xi) \delta(\xi-x) d\xi, \label{navier1} \\
\nabla \cdot u &= 0, \label{navier2} \\
u(\xi) &= \int u(x,t) \delta(x-\xi) dx =  u_B, \label{navier3}
\end{align}
where~$u$,~$p$ and~$f$ are the appropriately non-dimensionalized fluid velocity, pressure and surface force respectively. The force~$f$~acts as a Lagrange multiplier that imposes the no-slip boundary condition on the Lagrangian points~$\xi$, which arise from the discretization of a body moving with velocity~$u_B$. We consider the body to be a stationary flat plate at an angle of attack~$\alpha$; that is, here~$u_B = 0$. The Reynolds number is defined as $\Rey = Uc/\nu$ where $U$~is the free-stream velocity, $c$~is the chord-length and $\nu$~is the kinematic viscosity. Equations~(\ref{navier1}-\ref{navier3}) are discretized in space using a second-order finite-volume scheme on a staggered grid, and a discrete curl operation~$C^T(\cdot) \overset{\text{def}}{=} \nabla \times (\cdot)$ is introduced to eliminate the pressure and obtain a semi-discrete formulation in terms of the circulation~$\gamma$:
\begin{align}
\frac{d\gamma}{dt} + C^T E^T \tilde{f} &= -\beta C^T  C \gamma + C^T {\cal N} (q) + bc_{\gamma}, \label{semidisc} \\
ECs & = u_B = 0. \label{constraint}
\end{align}
The incompressibility condition~(\ref{navier2}) is implicitly satisfied by an appropriate construction of~$C$. The discrete Laplacian is represented by~$- C^T C \gamma$, using the identity~$\nabla^2 \gamma = \nabla (\nabla \cdot \gamma)- \nabla \times (\nabla \times \gamma) = - \nabla \times (\nabla \times \gamma)$; the constant $\beta = 1/\Rey \Delta^2$, where $\Delta$~is the uniform grid spacing. The operator~$E^T$ smears the Dirac delta function of~(\ref{navier1}) over a few grid points. The nonlinear term~${\cal N}(q)$ is the spatial discretization of~$q \times \gamma$, where $q$~is the discrete velocity flux, in turn related to the discrete stream function~$s$ and circulation~$\gamma$ as:
\begin{align}
q = C s, \qquad \gamma = C^T q, \qquad \mbox{and} \quad s = (C^T C)^{-1} \gamma. \label{stream}
\end{align}
%
%The incompressibility condition~(\ref{navier2}) is automatically satisfied as~$C$ is chosen to span the null-space of the discrete divergence operator; that is, if~$Dq = \nabla \cdot q$, then $DCs = 0$.
%
A uniform grid and a choice of simple boundary conditions result in a {\em fast} algorithm. With a uniform grid, the discrete Poisson equation~(\ref{stream}) is solved by means of the efficient discrete sine transform. The boundary conditions specified are Dirichlet and Neumann for the velocity components normal and tangential to the domain boundaries, which for the flow past a flat plate imply a uniform-flow in the far-field. These boundary conditions are however valid for only a sufficiently large domain, and are imposed by employing a computationally efficient multi-domain approach. The domain is considered to be embedded in a series of domains, each twice-as-large as the preceding, with a uniform but a {\em coarser} grid having the same number of grid points. The Poisson equation, with zero boundary conditions, is solved on the largest domain and the stream function is interpolated on the boundary of the smaller domain, which are in turn used to solve the Poisson equation on the smaller domain. For the flow past a flat plate considered here, the typical size of the largest domain is around 40~chord lengths in each direction. Finally, the time-integration  is performed using the implicit Crank-Nicolson scheme for the viscous terms and the second-order accurate, Adams-Bashforth scheme for the convective terms.

\subsection{Linearized and adjoint equations}  \label{sec:lin_adjoint}

For deriving reduced-order models useful for control design, we first linearize equation~(\ref{semidisc}) about a pre-computed steady state~($\gamma_0$,~$q_0$). The linearized equations are the same as equations~(\ref{semidisc},\ref{constraint}) with the nonlinear term~${\cal N}(q)$ replaced by its linearization about the steady state, and is denoted by ${\cal N}_L(\gamma_0) \gamma = q_0 \times \gamma + q \times \gamma_0$ where the flux~$q$ is related to~$\gamma$ by~(\ref{stream}). Thus, the linearized equations are:
\begin{align}
\frac{d\gamma}{dt} + C^T E^T \tilde{f} &= -\beta C^T  C \gamma + C^T {\cal N}_L(\gamma_0) \gamma,
\label{linear} \\
ECs & = 0.
\label{lin_constraint}
\end{align}
The boundary conditions for the linearized equations are~$bc_\gamma = 0$ on the outer boundary of the largest computational domain.

In order to derive the reduced-order models using the procedure described earlier, we need to perform adjoint simulations. In order to derive the adjoint equations, we define the following inner-product on the state-space:
\begin{align}
\ip<\gamma_1, \gamma_2> = \int_\Omega \gamma_1 \, (C^T C)^{-1} \gamma_2 \, dx.
\label{ip}
\end{align}
That is, the inner-product defined on the state-space is the standard $L^2$-inner product weighted with the inverse-Laplacian operator. This choice is convenient as it results in the adjoint equations which differ from the linearized equations only in the nonlinear term. A derivation is outlined in appendix~\ref{sec:app:adjoint} and the resulting equations are:
\begin{align}
\frac{d\zeta}{dt} + C^T E^T \psi & = - \beta C^T C \zeta + (C^T C){\cal N}_L(\gamma_0)^T q_a , \label{adj1} \\
EC \xi & = 0, \label{adj2}
\end{align}
where the variables~$\zeta$,~$\xi$ and~$\psi$ are the duals of the discrete circulation~$\gamma$, stream function~$s$ and body force~$\tilde{f}$, respectively, and~$q_a = C \xi$ is the dual of flux~$q$. The adjoint of the linearized nonlinear term is~$(C^T C){\cal N}_L(\gamma_0)^T q_a$, which can be shown to be a spatial discretization of~$\nabla \times (\gamma_0 \times q_a) - \nabla^2(q_0 \times q_a)$. Since equation~(\ref{adj1}) differs from~(\ref{linear}) only in the last term on the right hand side, the numerical integrator for the adjoint equations can be obtained by a small modification to the linearized equations solver.

The nature of the multi-domain scheme used to approximate the boundary conditions of the smallest computational domain, results in a multi-domain discrete Laplacian which is not exactly self-adjoint to numerical precision. As a result, the adjoint formulation given by~(\ref{adj1},\ref{adj2}) which also uses the same multi-domain approach, is not precise and results in small, rather insignificant, errors in the computation of the reduced-order models.

\begin{figure}
\centering
\includegraphics[height=2in]{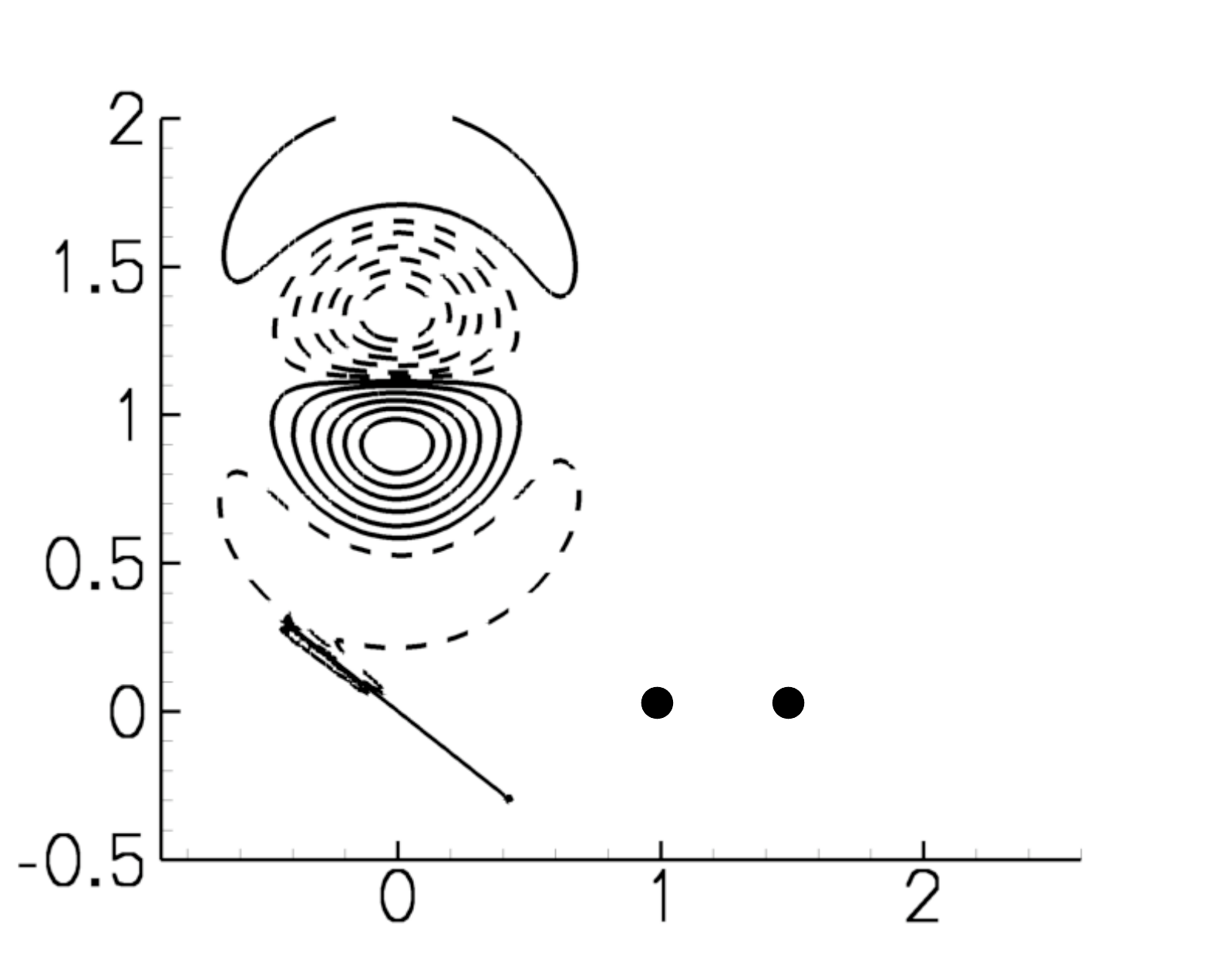}
\caption{ Actuator modeled as a localized body force near the leading edge of the flat plate, with the angle of attack fixed at~$\alpha = 35^\circ$. Vorticity contours are plotted, with negative contours shown by dashed lines. The velocity-sensor locations are marked by solid circles.}
\label{fig:actuator_sensor}
\end{figure}

\section{Results: flow past a flat plate} \label{section:model_results}

We apply the model reduction techniques developed in the previous sections to the uniform flow past a flat plate in two spatial dimensions, at a low Reynolds number, $\Rey = 100$. We obtain reduced-order models of a system actuated by means of a localized body force near the leading edge of the flat plate; the vorticity contours of the flow field obtained on an impulsive input to the actuator are shown in Fig.~\ref{fig:actuator_sensor}. Using these reduced-order models, we develop feedback controllers that stabilize the unstable steady state at high angles of attack. We first assume full-state feedback, but use output projection described in section~\ref{sec:outproj} to considerably decrease the number of outputs in order to make the model computation tractable. Later, we relax the full-state feedback assumption, and develop a more practical observer-based controller which uses a few velocity measurements in the near-wake of the flat plate~(shown in Fig.~\ref{fig:actuator_sensor}) to reconstruct the entire flow.

\subsection{Numerical parameters} \label{sec:params}
The grid size used is $250\times 250$, with the smallest computational domain given by~$[-2,3]\times[-2.5,2.5]$, where lengths are non-dimensionalized by the chord of the flat plate, with its center located at the origin.  We use 5 domains in the multiple-grid scheme, resulting in an effective computational domain~$2^4$ times larger the size of the smallest domain; thus the largest domain is given by~$[-32,48]\times[-40,40]$. The timestep used for all simulations is~$\delta_t = 0.01$.

\subsection{Steady-state analysis} \label{sec:steady_states}
Since our approach is to obtain reduced-order models of the flow linearized about a given steady state, we first need to compute these steady states. The model-reduction of unstable systems involves projecting the dynamics onto a stable subspace, for which we also need to compute the right and left eigenvectors of the linearized dynamics. This section concerns this steady-state analysis, using a ``timestepper-based'' approach as outlined in~\cite{TuBar99} and~\cite{kelley04}.

A simple way of computing stable steady states is by simply evolving the time-accurate simulation to stationarity. However, unstable steady states cannot be found in this manner, and stable steady states near a bifurcation point could take very long to converge. Instead, we use a timestepper-based approach which involves writing a computational wrapper around the original computational routine to compute the steady states using a Newton iteration. If the numerical timestepper advances a circulation field~$\gamma^k$ at a timestep~$k$ to a circulation field $\gamma^{k+T}\equiv \Phi_T(\gamma^k)$ after $T$~timesteps, the steady state is given by the field~$\gamma_0$ that satisfies
\begin{equation}
g(\gamma_0) =  \gamma_0 - \Phi_T(\gamma_0)=0.
\label{newton}
\end{equation}
The steady states are given by zeros of~$g(\gamma_0)$, which could, in principle, be solved for using Newton's method. However, the standard Newton's method involves computing and inverting Jacobian matrices at each iteration, which is computationally infeasible due to the large dimension of fluid systems. Instead of computing the Jacobian, we use a Krylov-space based iterative solver called Generalized Minimal Residual Method~(GMRES) developed by~\cite{Saad-86} to compute the Newton update (see~\cite{Kelley-95} and~\cite{Trefethen-97} for a description of the method). This method requires computation of only Jacobian-vector products~$Dg(\gamma)\cdot v$, which can be approximated using finite differences as~$[g(\gamma + \epsilon v) - g(\gamma)]/\epsilon$, for~$0<\epsilon \ll 1$. So, the Jacobian-vector products can also be computed by invoking the appropriately-initialized timestepper.
A nice feature of GMRES is relatively fast convergence to the steady state when the eigenvalues of the Jacobian~$Dg(\gamma_0)$ occur in clusters; see \cite{Kelley-95} and \cite{kelley04} for details. For systems with multiple time-scales, such as Navier-Stokes, most of the eigenvalues of the continuous Jacobian lie in the far-left-half of the complex plane. Thus, the corresponding eigenvalues of the discrete Jacobian~$D\Phi_T$, for a sufficiently large value of~$T$, cluster near the origin.

\begin{figure}
\centering
\includegraphics[height=1.4in]{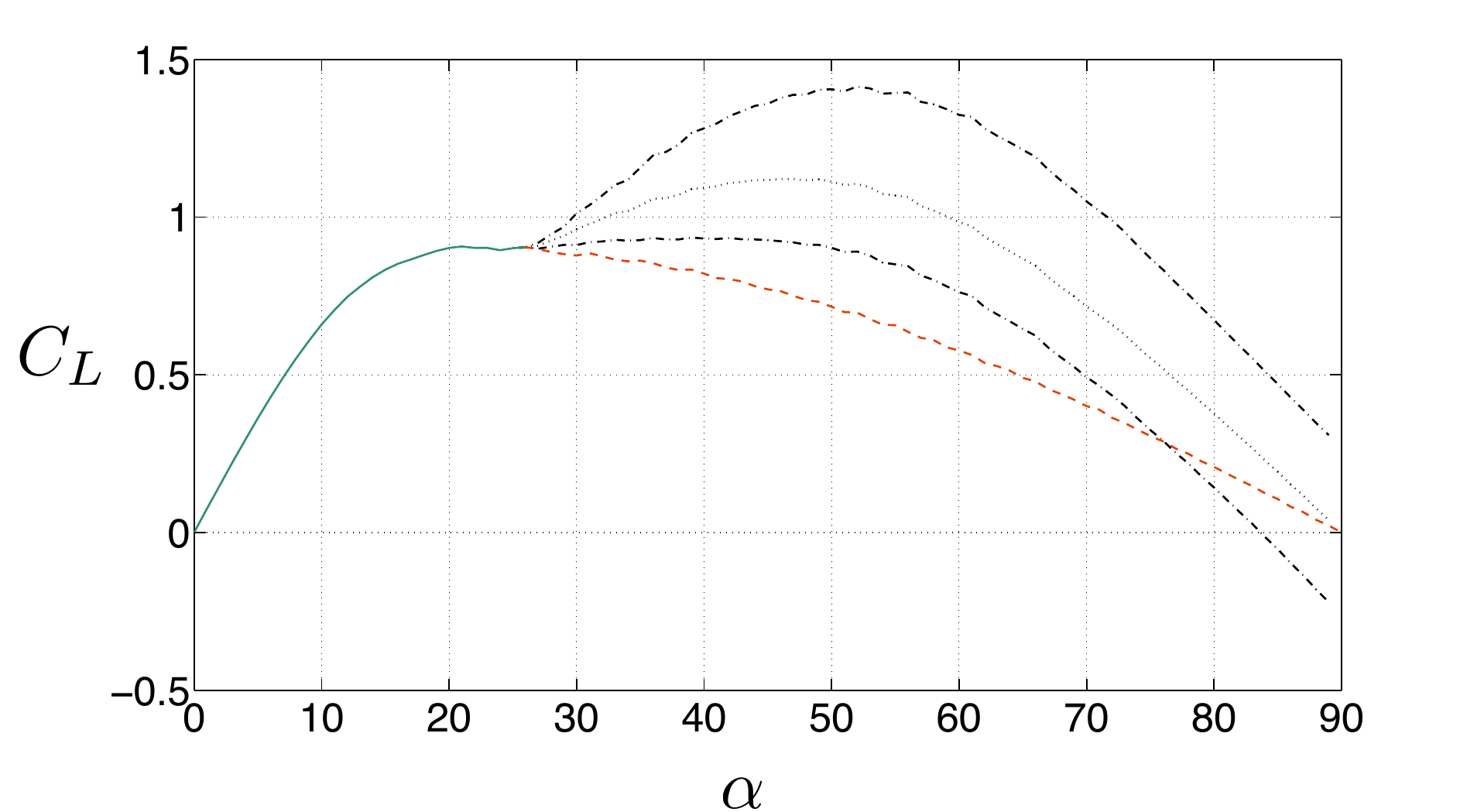}
\includegraphics[height=1.4in]{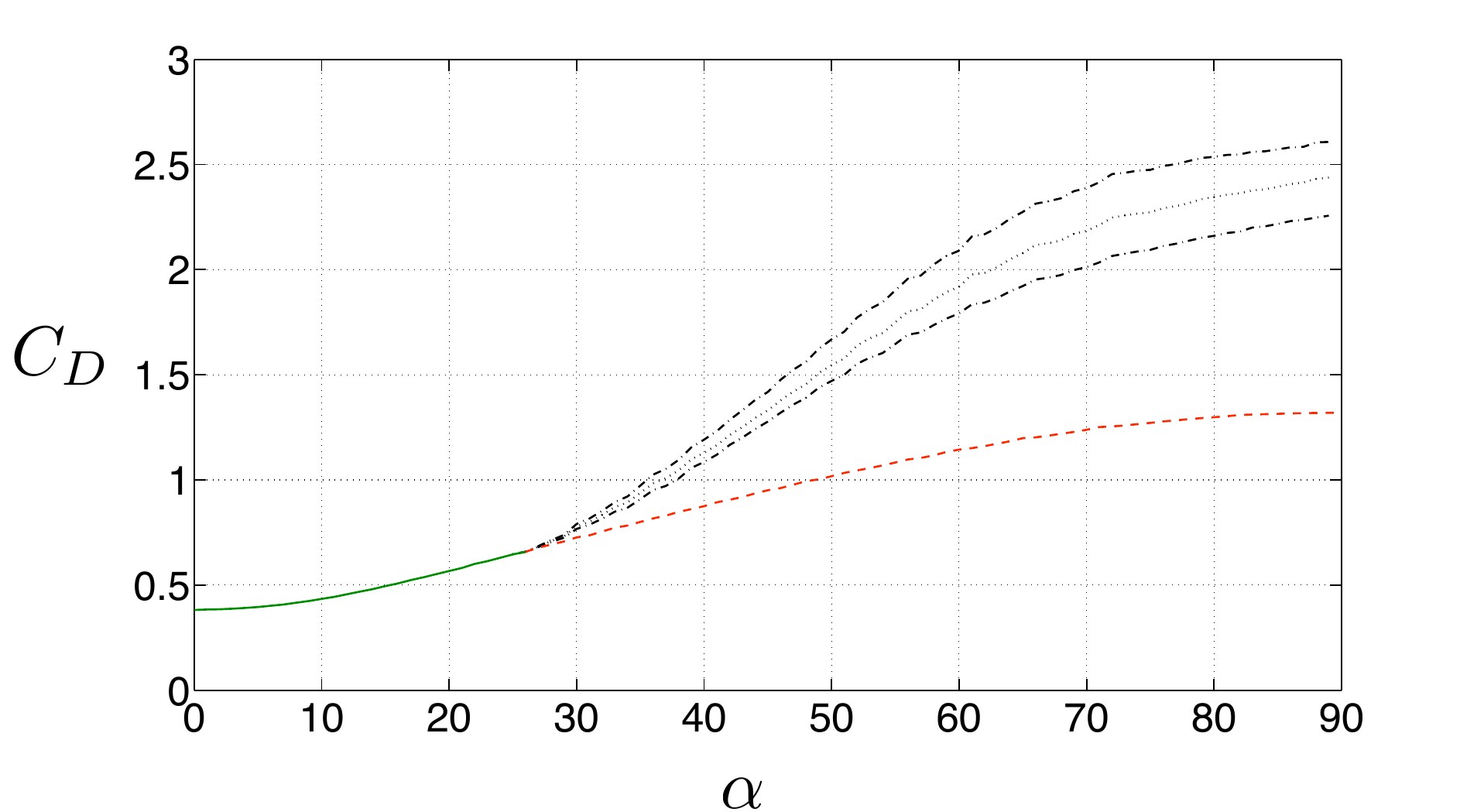}
\includegraphics[height=1.9in]{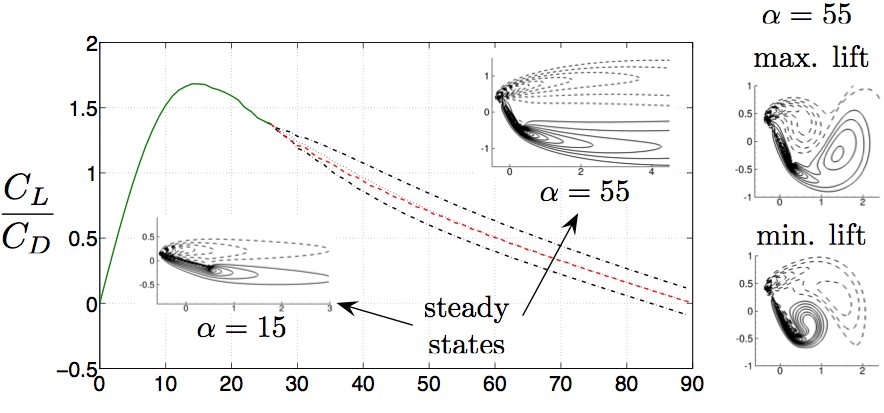}
\caption{Forces on a flat plate at a fixed angle of attack~$\alpha$ and at~$\Rey = 100$, showing a transition from a stable equilibrium to periodic vortex shedding at~$\alpha \approx 26$. Shown are the force coefficients corresponding to the stable~({\color{green}$\solid$}) and unstable~({\color{red}$\dashed$}) steady states, and the maximum and minimum~($\chndot$), and the mean~($\dotted$) values during periodic vortex shedding. Also shown are the vorticity contours (negative values in dashed lines) of steady states at~$\alpha = 15^\circ, 55^\circ$ and the flow fields corresponding to the maximum and minimum force coefficients at~$\alpha = 55^\circ$. }
\label{fig:bifurc}
\end{figure}

\begin{figure}
\centering \hspace{0.5in} (a)  \hspace{2.5in}  (b) \\
\includegraphics[height=1.5in]{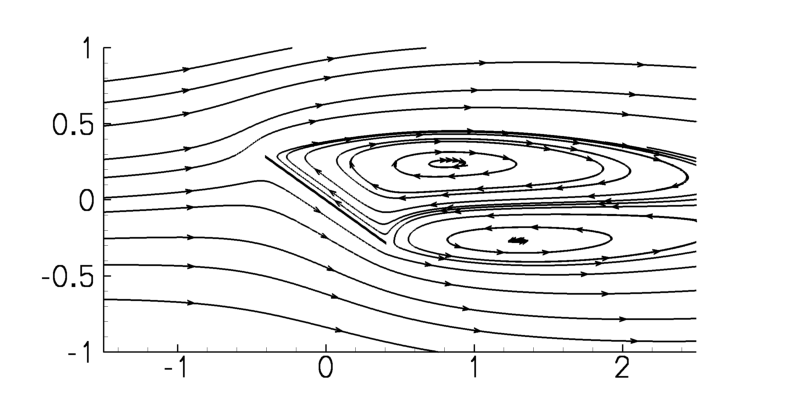}
\includegraphics[height=1.5in]{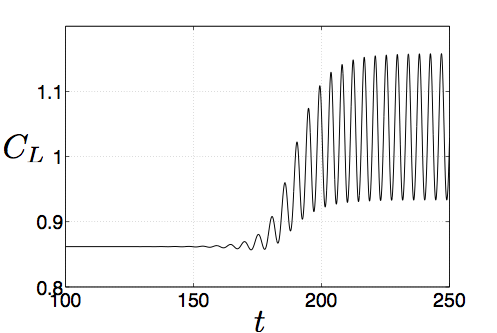}
\caption{(a) Streamlines of the unstable steady state at $\alpha = 35$. (b) $C_L$ vs.~time, with the steady state as an initial condition.} \label{fig:ss}
\end{figure}

The procedure described above is used to compute the branch of steady states for the angles of attack~$0<\alpha < 90^\circ$; the parameter~$T$  in~(\ref{newton}) is fixed to 50~timesteps. The lift and drag coefficients,~$C_L$ and~$C_D$, and their ratio~$C_L/C_D$ with changing~$\alpha$ are plotted in Fig.~\ref{fig:bifurc}. As with flow past bluff bodies with increasing Reynolds number (for example, see~\cite{PrMaBoyer-jfm87}), the flow undergoes a Hopf bifurcation from a steady flow to periodic vortex shedding as the angle of attack~$\alpha$ is increased beyond a critical value~$\alpha_c$, which in our computations is~$\alpha_c \approx 27^\circ$. Also plotted in the figure are the maximum, minimum, and mean values of the forces during shedding for~$\alpha > \alpha_c$. We see that the (unstable) steady state values of the lift coefficient are smaller than the minimum for the periodic shedding till~$\alpha \approx 75^\circ$, after which they are slightly higher, but still smaller than the mean lift for the periodic shedding.  The (unstable) steady state drag is much lower than the minimum value for periodic shedding. The ratio~$C_L/C_D$ of the (unstable) steady state is close to the mean value for shedding. Thus, if the large fluctuations in the forces are undesirable at high angles of attack, it would be useful to stabilize the unstable state. The steady state at~$\alpha = 35^\circ$ is shown in Fig.~\ref{fig:ss}(a), and a time history of the lift coefficient~$C_L$ with this steady state as an initial condition is shown in Fig.~\ref{fig:ss}(b). Since the steady state is unstable, the numerical perturbations excite the instability, and the flow eventually transitions to periodic vortex shedding.

\begin{figure}
\centering
\includegraphics[height=1.2in]{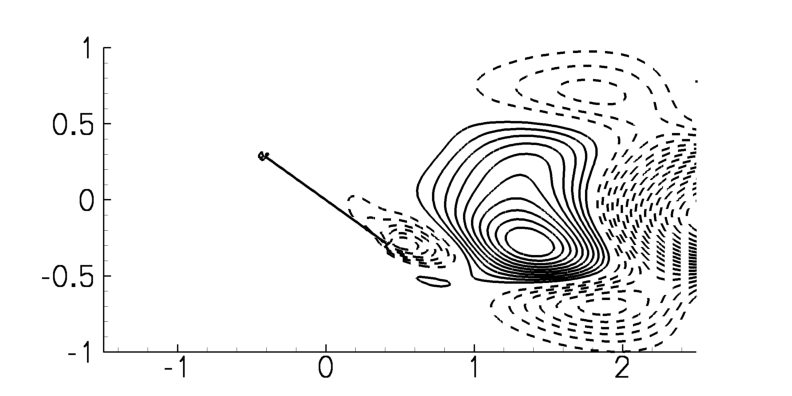}
\includegraphics[height=1.2in]{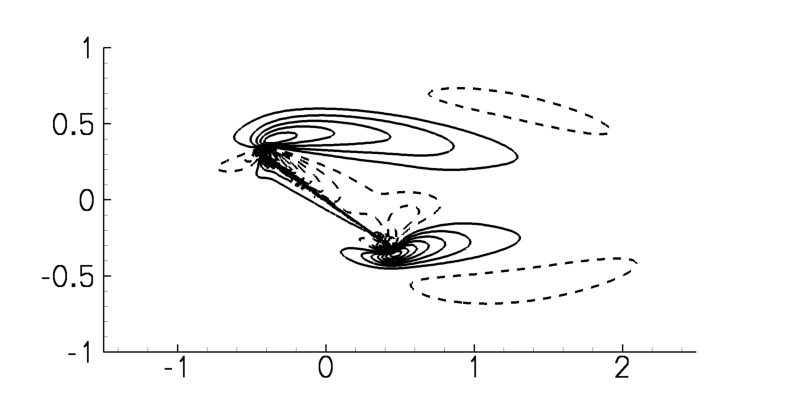} \\
\includegraphics[height=1.2in]{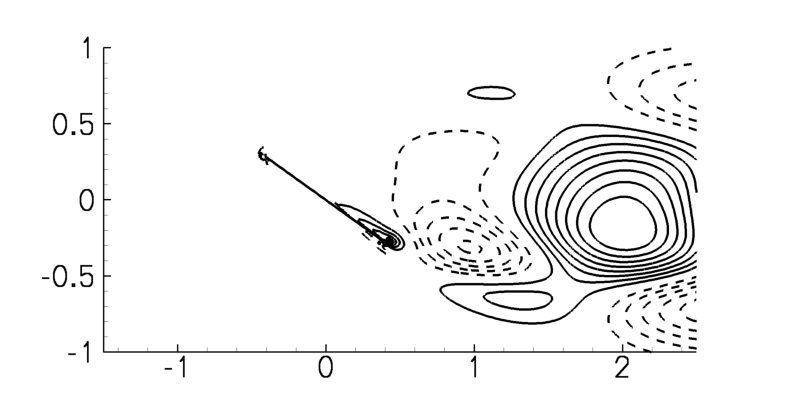}
\includegraphics[height=1.2in]{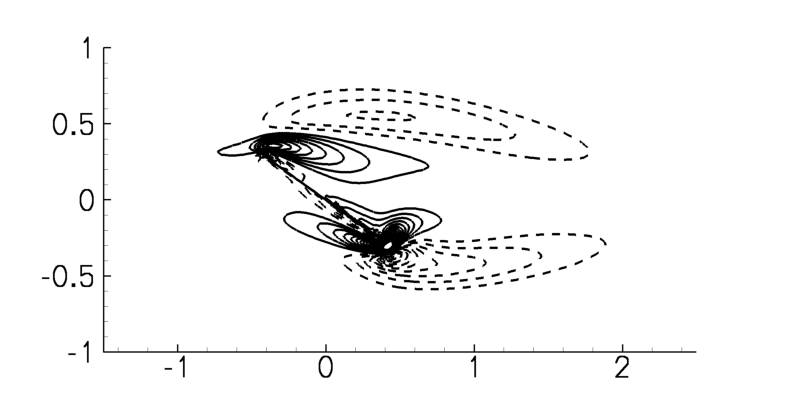}
\caption{Basis vectors of the unstable eigenspace of the linearized (left) and the adjoint (right) equations. Vorticity contours are plotted (negative contours are dashed).}
\label{fig:evec}
\end{figure}

We also compute a basis spanning the right and left unstable eigenspaces~($\Phi_u$ and~$\Psi_u$) of the flow linearized about the unstable steady states, which are required in our model reduction procedure, for restricting dynamics onto the stable subspace. As the flow undergoes a Hopf bifurcation, a complex pair of eigenvalues crosses the imaginary axis from the left half of the complex plane; thus the dimension of the unstable subspace is two. Here, we obtain the basis spanning this subspace by a numerical implementation of the power method~(see page~191 of~\cite{Trefethen-97}). We begin the linearized simulation with a very small random noise to excite the instability. After a sufficiently long time, the stable modes decay, and the dynamics lies close to the unstable subspace. Any two independent snapshots of the remaining dynamics gives a good approximation of the basis spanning the unstable eigenspace. Similarly, a basis spanning the left unstable eigenspace can be computed by initializing the adjoint equations with a small random noise. A basis spanning the right and left unstable eigenspaces of the flow linearized about the steady state at~$\alpha = 35^\circ$ is plotted in Fig.~\ref{fig:evec}. These modes are qualitatively similar to the structures during periodic vortex shedding, but have different spatial wavelengths, as reported in earlier studies by \cite{Noack-03} and~\cite{Barkley-06}.

\begin{figure}
\centering
(a) \\
\includegraphics[height=1.3in]{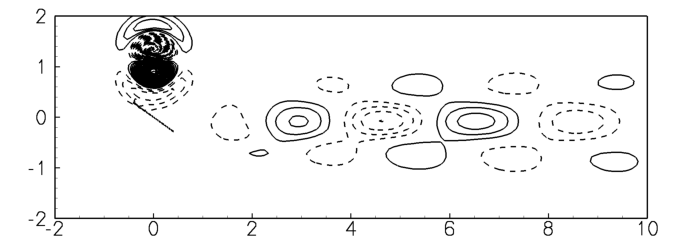} \\
(b) \hspace{2.5in} (c) \\
\includegraphics[height=1.3in]{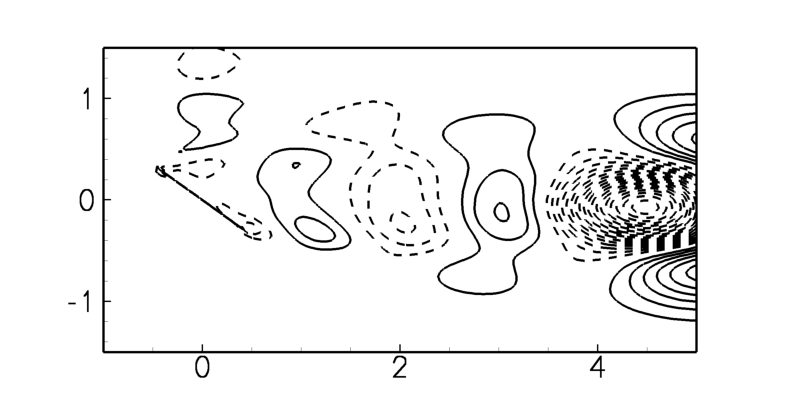}
\includegraphics[height=1.3in]{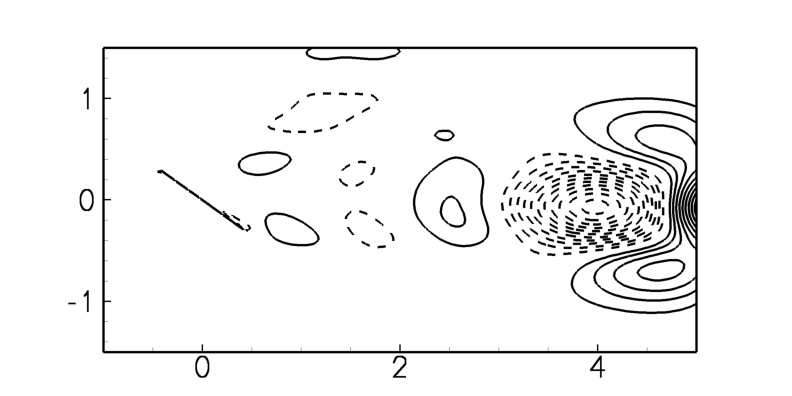}
\caption{Vorticity contours of (a) the flow field shown in Fig.~\ref{fig:actuator_sensor}, projected onto the stable subspace, and (b,c) the first- and fifth-most energetic POD modes of the impulse response, restricted to the stable subspace.}
\label{fig:actuator}
\end{figure}

\subsection{Reduced-order models} \label{sec:models}

We now describe the process involved in deriving reduced-order models of the input-output response of~(\ref{ss}), which in this example are the linearized incompressible Navier-Stokes equations~(\ref{linear}, \ref{lin_constraint}). The actuator used is a localized body force close to the leading edge of the flat plate, plotted in Fig.~\ref{fig:actuator_sensor}. The models are derived using the procedure outlined in section~\ref{sec:approx_baltrunc_unstable}. As seen in equation~(\ref{rom}), the output of the system is considered to be the entire velocity field, observed as a projection onto (a) the unstable eigenspace, and (b) the span of the leading POD modes of the impulse response restricted to the stable subspace.

The first step in computing the reduced-order models is to project the flow field~$B$ onto the stable subspace of~(\ref{linear}, \ref{lin_constraint}) using the projection operator~$\mathbb{P}_s$ defined in equation~(\ref{proj_stable}); the unstable eigenvectors computed in section~\ref{sec:steady_states} are used to define~$\mathbb{P}_s$ numerically. The vorticity contours of the corresponding flow field~$\mathbb{P}_s B$ are plotted in Fig.~\ref{fig:actuator}a. The next step is to compute the impulse response of~(\ref{ss_stable}). Instead, for practical reasons, we compute the impulse response of
\begin{align}
\dot{x_s} = \mathbb{P}_s A  x_s + \mathbb{P}_s B u;
\label{ss_stable_numerical}
\end{align}
that is, at each timestep of integration, we project the state~$x_s$ onto the stable subspace of~$A$. Because the stable subspace is an invariant subspace for the linearized dynamics~(\ref{linear}), theoretically, the impulse responses of equations~(\ref{ss_stable}) and~(\ref{ss_stable_numerical}) are exactly the same, and they are the same as that obtained by restricting the impulse response of~(\ref{ss}) to its stable subspace. However, due to the (small)~numerical inaccuracy of the projection~$\mathbb{P}_s$ (which is a result of the numerical inaccuracy of the unstable eigenspaces~$\Phi_u$ and~$\Psi_u$), the dynamics of~(\ref{ss_stable}) is not strictly restricted to the stable subspace and, in the long term, grows without bound in the unstable direction. Next, we compute the POD modes~$\theta_s^i$ of the impulse response of~($\ref{ss_stable_numerical}$), and consider the output of~(\ref{ss_stable_numerical}) to be the state~$x_s$ projected onto a certain number of these POD modes. Here, 200 snapshots spaced every 50~timesteps were used to compute the POD modes. The leading 4 and 10 POD modes contain~$84.47\%$ and~$99.98\%$ of the energy respectively and, as it has been observed in previous studies (see \cite{DeKeKaOrszag-91} and \cite{IlakRowley-pof08}), these modes come in pairs in terms of their energy content, a characteristic of traveling structures; the leading first and third POD modes are shown in Fig.~\ref{fig:actuator}.

\begin{figure}
\centering
\includegraphics[height=1.3in]{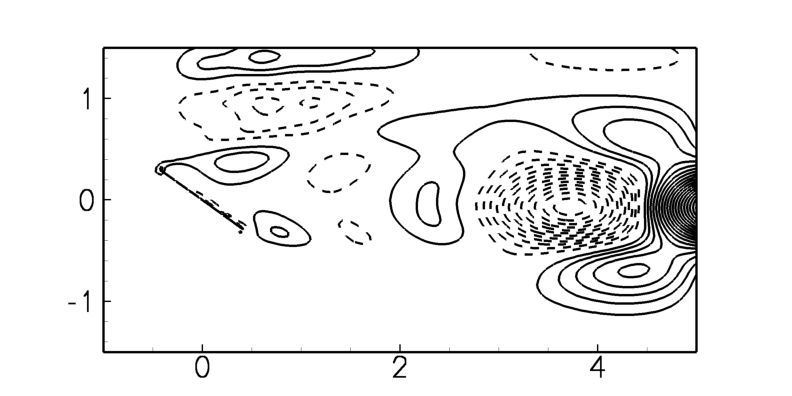}
\includegraphics[height=1.3in]{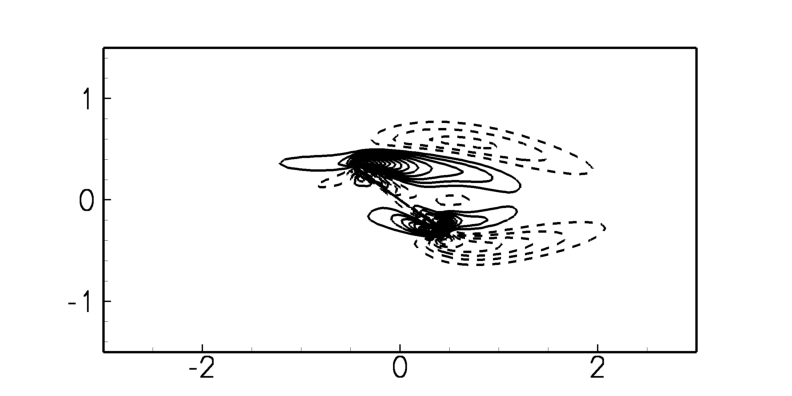} \\
\includegraphics[height=1.3in]{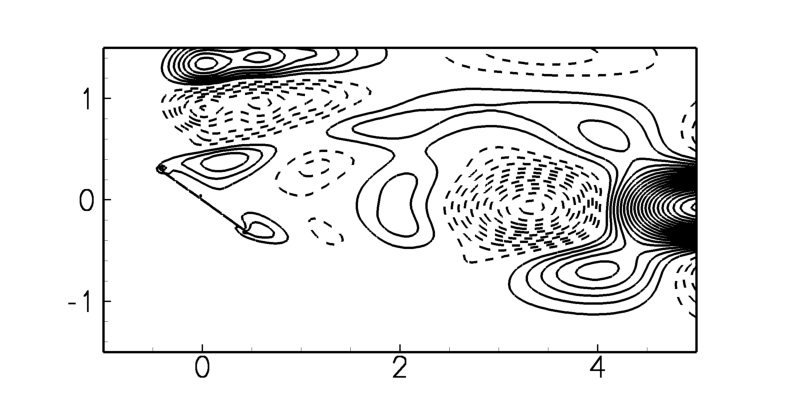}
\includegraphics[height=1.3in]{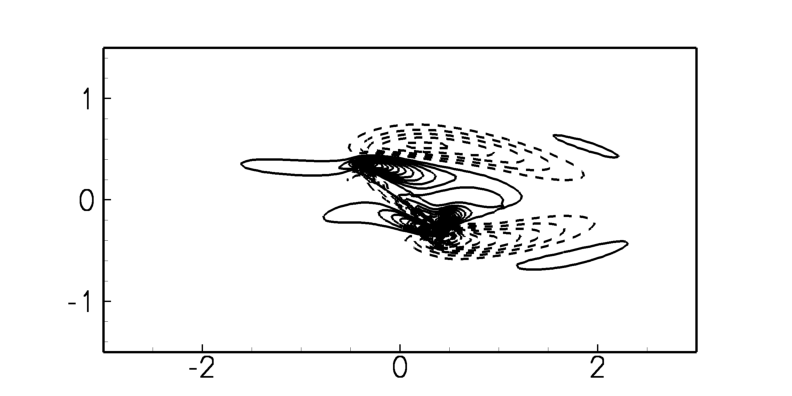}
\caption{Vorticity contours of the leading (in the order of Hankel singular values of the stable subspace dynamics) first and third balancing (left) and adjoint (right) modes.}
\label{fig:modes}
\end{figure}

The next step is to compute the adjoint snapshots, with the POD modes of the impulse response (projected onto the stable subspace of the adjoint) as the initial conditions. As the linearized impulse response, these simulations are also restricted to the stable subspace. Again, instead of computing the response of~(\ref{ss_adj_stable}), we compute that of the following system:
\begin{align}
\dot{z_s} = \mathbb{P}_s^\ast A^\ast + \mathbb{P}_s^\ast C^\ast v.
\label{ss_adj_stable_numerical}
\end{align}
The snapshots of the impulse responses of systems~(\ref{ss_stable_numerical}) and~(\ref{ss_adj_stable_numerical}) are stacked as columns of~$X$ and~$Z$, and using the expressions~(\ref{svd}) and~(\ref{bal_transf}), we obtain the balancing modes~$\phi_s^i$ and the adjoint modes~$\psi_s^i$.
We used 200~snapshots of the linearized simulation and 200 snapshots of each adjoint simulation, with the spacing between snapshots fixed to 50~timesteps, to compute the balancing transformation. These number of snapshots and the spacing were sufficient to accurately compute the modes; further reduction in the spacing did not significantly change the singular values from the SVD computation~(\ref{svd}). We considered the outputs to be a projection onto 4, 10 and 20~POD modes (corresponding to 4, 10 and 20~mode {\em output-projections}, as introduced in section~\ref{sec:outproj}). Using these modes, we use the expressions in equation~(\ref{rom_compact},\ref{output_compact2}) to obtain the matrices~$\widetilde{A}_s, \widetilde{B}_s, \widehat{C}_s$ defining the reduced-order model of the stable-subspace dynamics. The vorticity contours of the balancing and the adjoint modes, for a 10-mode output projected system, are plotted in~Fig.~\ref{fig:modes}.
%The comparison of two subspaces can be quantified as well, and it has been observed in previous studies that the POD and balancing modes span very similar subspaces.\cite{IlakRow07}
The adjoint modes provide a direction for projecting the linearized equations onto the subspace spanned by the balancing modes. Since these modes are quite different from the POD and the balancing modes, the resulting models are also quite different from those obtained using the standard POD-Galerkin technique wherein an orthogonal projection is used. Since the models obtained using balanced truncation are known to perform better than the POD-Galerkin models, as reported by \cite{IlakRowley-pof08}, the better performance could be attributed to a better choice of projection using the adjoint modes.

\begin{figure}
\centering
\includegraphics[height=2.4in]{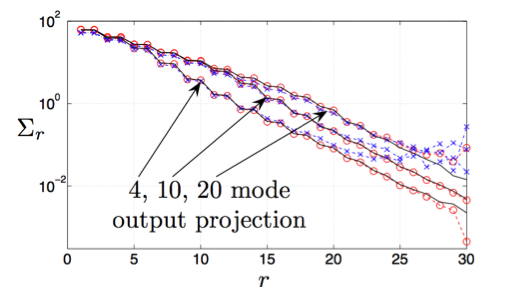}
\caption{The empirical Hankel singular values~($\solid$) and the diagonal elements of the controllability~({\color{red}$\dashed$,~$\circ$}) and observability~({\color{blue}$\chndash$,~$\times$}) Gramians of  a 25-mode model with a 4, 10, and 20-mode output projection, for the unstable steady state at~$\alpha=35$.}
\label{fig:unstable_hsv}
\end{figure}

Since the reduced-order models of the stable-subspace dynamics are approximately balanced, the controllability and observability Gramians of the $a_s$-dynamics of~(\ref{rom_compact}), given by expressions~(\ref{gramians_stable}), are approximately equal and diagonal. Further, their diagonal values are approximately the same as the Hankel singular values~$\sigma_i$ obtained by the SVD~(\ref{svd}). The diagonal values of the Gramians and the singular values for different output projections are plotted in Fig.~\ref{fig:unstable_hsv} for a 30-state reduced-order model.  With increasing order of output projection, the HSVs converge to the case with full-state output, and the number of converged HSVs is roughly equal to the order of output projection, as was observed by~\cite{IlakRowley-pof08}. We see that the diagonal elements of both the Gramians are very close to the HSVs for the first 20 modes. For higher modes, the observability Gramians are inaccurate, which is due to a small inaccuracy of the adjoint formulation mentioned in section~\ref{sec:lin_adjoint}. For controller design, we use models of order~$\leq 20$, for which these Gramians are sufficiently accurate.

\begin{figure}
\centering
\includegraphics[height=1.7in]{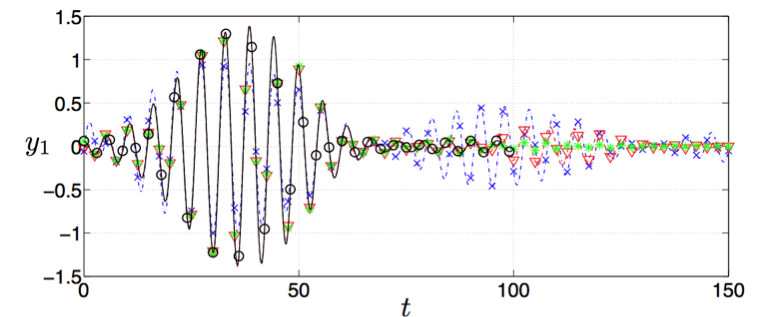}
\includegraphics[height=1.7in]{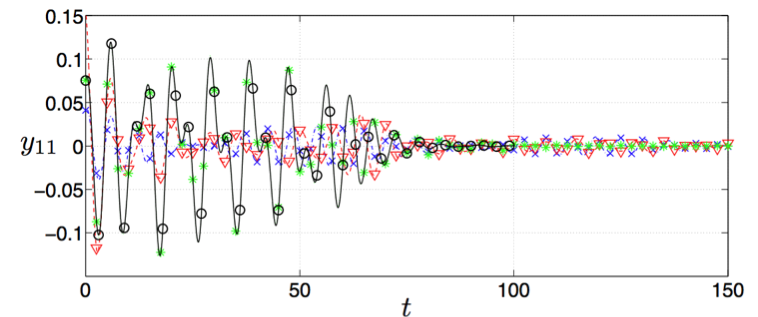}
\caption{Outputs (projection of the flow field onto POD modes) from a reduced-order model obtained using a 20-mode output projection. The first (top figure) and eleventh (bottom figure) outputs of the DNS~($\solid$, $\circ$) are compared with predictions of models with 4~({\color{blue}$\chndash$, $\times$}), 10~({\color{red}$\dashed$, $\triangledown$}), and 20~({\color{green}$\dotted$, $\ast$})~modes.}
\label{fig:outputs_stable_oproj10}
\end{figure}

Finally, to test the accuracy of the reduced-order models, we compare the impulse responses of system~(\ref{ss_stable_numerical}) (that is, restricted to the stable subspace) with that of the model~(\ref{rom_compact}), restricting~$a_u = 0$. In particular, we compare the outputs of the two systems, which are the projection onto the POD modes; a representative case in Fig.~\ref{fig:outputs_stable_oproj10} shows the results of 4, 10 and 20 mode models of a system approximated using a 20-mode output projection (the outputs are projection onto the leading 10~POD modes). The first output, which is a projection onto the first POD mode, is well captured by all the models until~$t \approx 60$, while the 20-mode model performs well for all time. Also shown is the eleventh output, which is well captured only by the 20-mode model. As we will see later, it is important to capture the higher-order outputs for design of observers.

\begin{figure}
\centering
\includegraphics[height=1.3in]{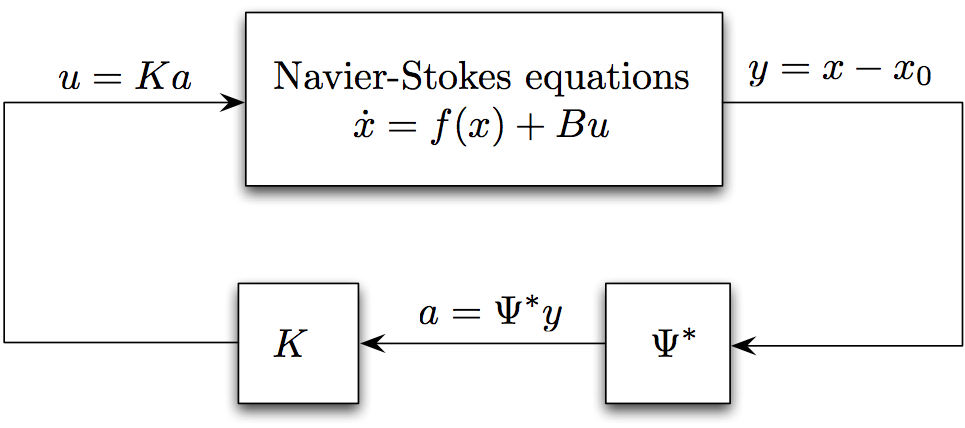}
\caption{Schematic of the implementation of full-state feedback control in the nonlinear simulations. The entire velocity is first projected onto the unstable eigenvectors and the stable subspace POD modes to compute the reduced-order state~$a$. The state is then multiplied by the gain~$K$, computed based on the reduced-order model using LQR, to obtain the control input~$u$.}
\label{fig:cartoon_fullstate_control}
\end{figure}

\begin{figure}
\centering
%{\includegraphics[height=1.4in]{impulse_outproj4_ustb1}
%\includegraphics[height=1.4in]{impulse_outproj4_stb3}}
\includegraphics[height=1.4in]{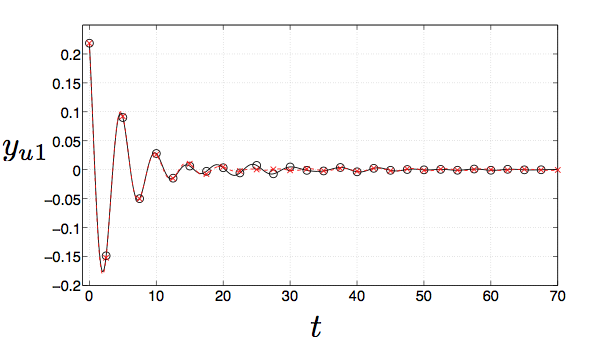}
\includegraphics[height=1.4in]{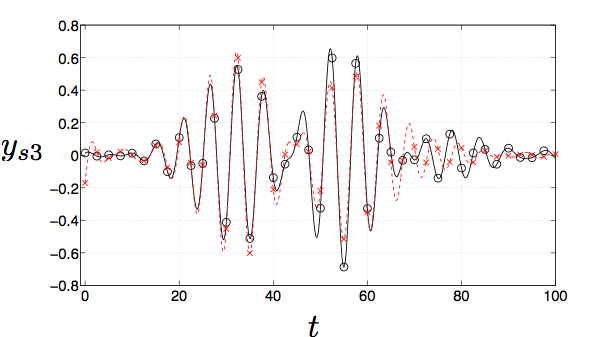}
\caption{Control gain obtained using LQR, and the initial condition is that obtained by an impulsive input to the system. Control is turned on at~$t=0$. Comparison of the outputs~$y_{u1}$ and~$y_{s3}$ of a 12-mode reduced-order model ({\color{red}$\dashed,\times$}) with the projection of data from the linearized simulation ($\solid,\circ$).} \label{fig:contlin}
\end{figure}

\subsection{Full-state feedback control} \label{sec:fullstatefb}

The resulting models can now be used along with standard linear control techniques to obtain stabilizing controllers. We use Linear Quadratic Regulator~(LQR) to compute the gain~$K$ so that the eigenvalues of~$(\widetilde{A}+\widetilde{B}K)$ (where the matrices were defined in~(\ref{rom_compact})) are in the left-half of the complex plane, and the input~$u=Ka$ minimizes the cost
\begin{align}
J[a,u] = \int_0^\infty \, (a^\ast Q a + u^\ast R u) \, dt,  \label{lqr_cost}
\end{align}
where~$Q$ and~$R$ are positive weights computed as follows. We choose~$Q$ such that the first term in the integrand of~(\ref{lqr_cost}) represents the energy, that is, we use~$Q=\widetilde{C}^\ast \widetilde{C}$, with~$\widetilde{C}$ defined in~(\ref{output_compact}). The weight~$R$ is chosen to be a multiple of the identity~$cI$, and typically~$c$ is chosen to be a large number~$\sim 10^{4-7}$, to avoid excessively aggressive controllers. The control implementation steps are sketched in Fig.~\ref{fig:cartoon_fullstate_control}; first compute the reduced-order state~$a$, using the expression~(\ref{init_model}), then the control input is given by~$u = Ka$. Here, we derive the gain~$K$ based on a 12-mode reduced-order model (with 2 unstable and 10 stable modes), using~$R=10^{5}$, and include the same in the original linearized and nonlinear simulations. The output is approximated using a 4-mode output projection. The difference between the linear and nonlinear simulations is that, in the latter, the steady state field~$x_0$ is subtracted from the state~$x$, before projecting onto the modes to compute the reduced-order state~$a$.

Fig.~\ref{fig:contlin} compares the model predictions with the projection of data from the simulations of the linearized system~(\ref{linear}, \ref{lin_constraint}), with a control input. The initial condition used is the flow field obtained from an impulsive input to the actuator. Both the states shown in the figure eventually decay to zero, which implies that the perturbations decay to zero, thus stabilizing the unstable steady state. More importantly, the model predicts the outputs accurately for the time horizon shown in the plots.

\begin{figure}
\centering
%{\includegraphics[width=4.5in]{control_oproj4_lift}}
{\includegraphics[width=4.5in]{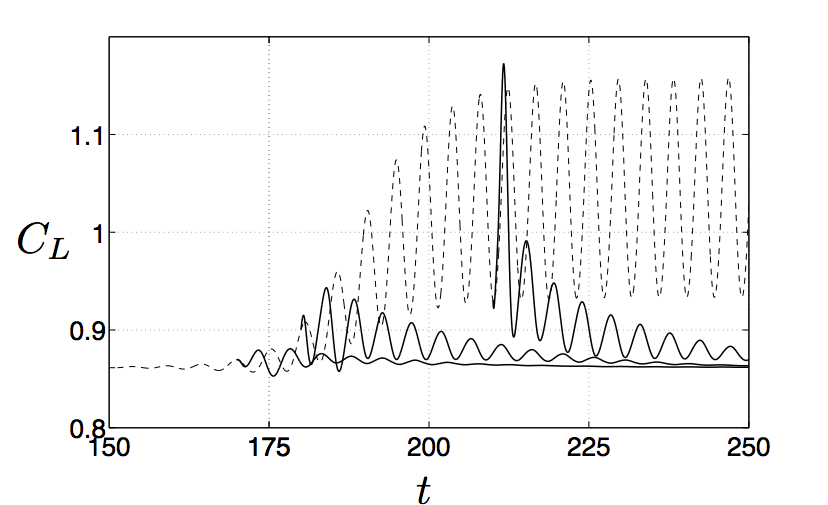}}
\caption{Lift-coefficient~$C_L$ vs. time~$t$, for full-state feedback control, with control turned on at different times in the base uncontrolled simulation. The base case with no control ($\dashed$) has the unstable steady state as the initial condition, and transitions to periodic vortex shedding. The control is tested for different initial conditions, corresponding to~$t = 170, 180, 210$ of the base case, and stabilizes the steady state in all the cases ($\solid$).} \label{fig:cont_fullfback_lift}
\end{figure}

\begin{figure}
\centering
\includegraphics[height=1.4in]{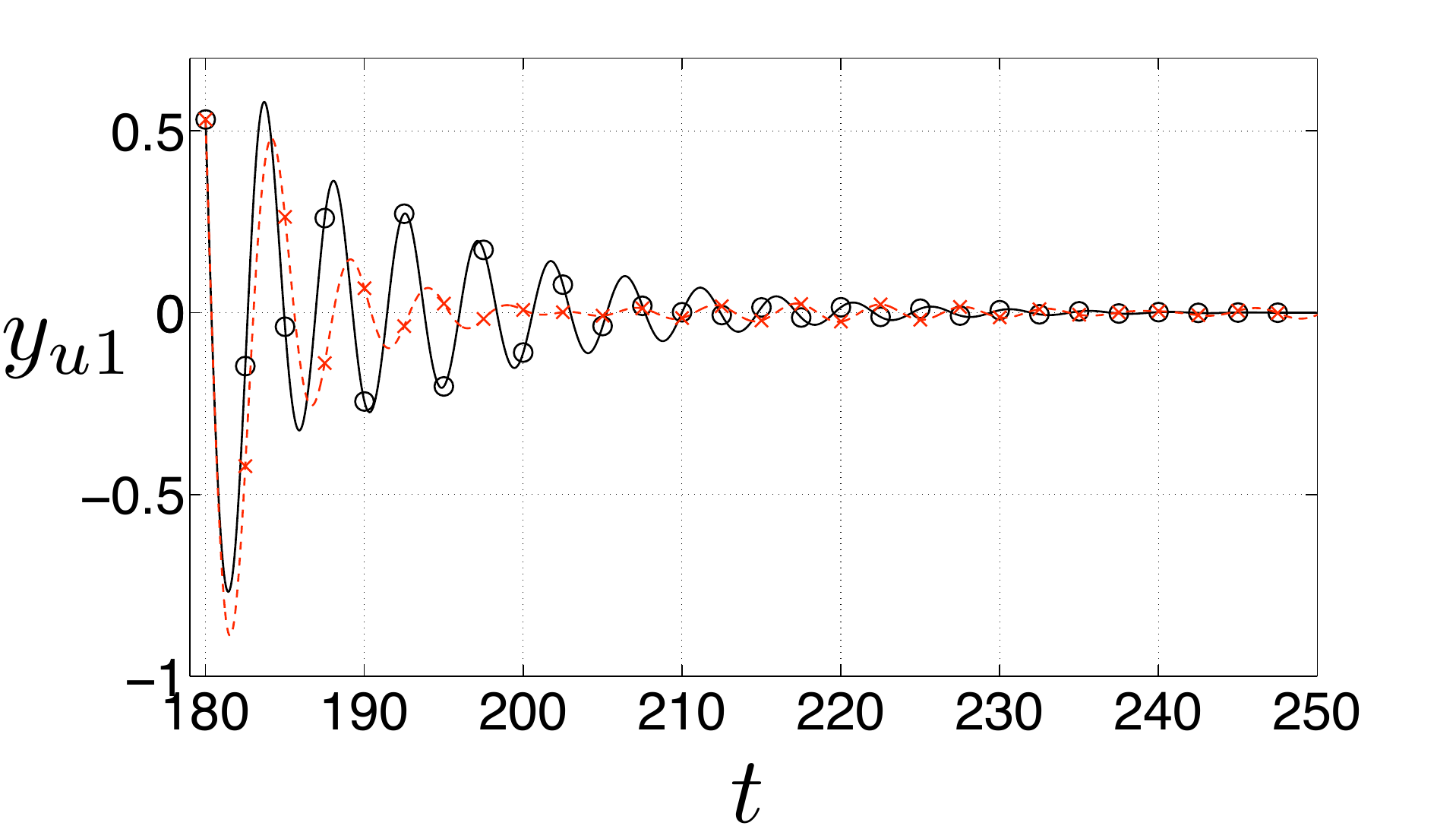}
\includegraphics[height=1.4in]{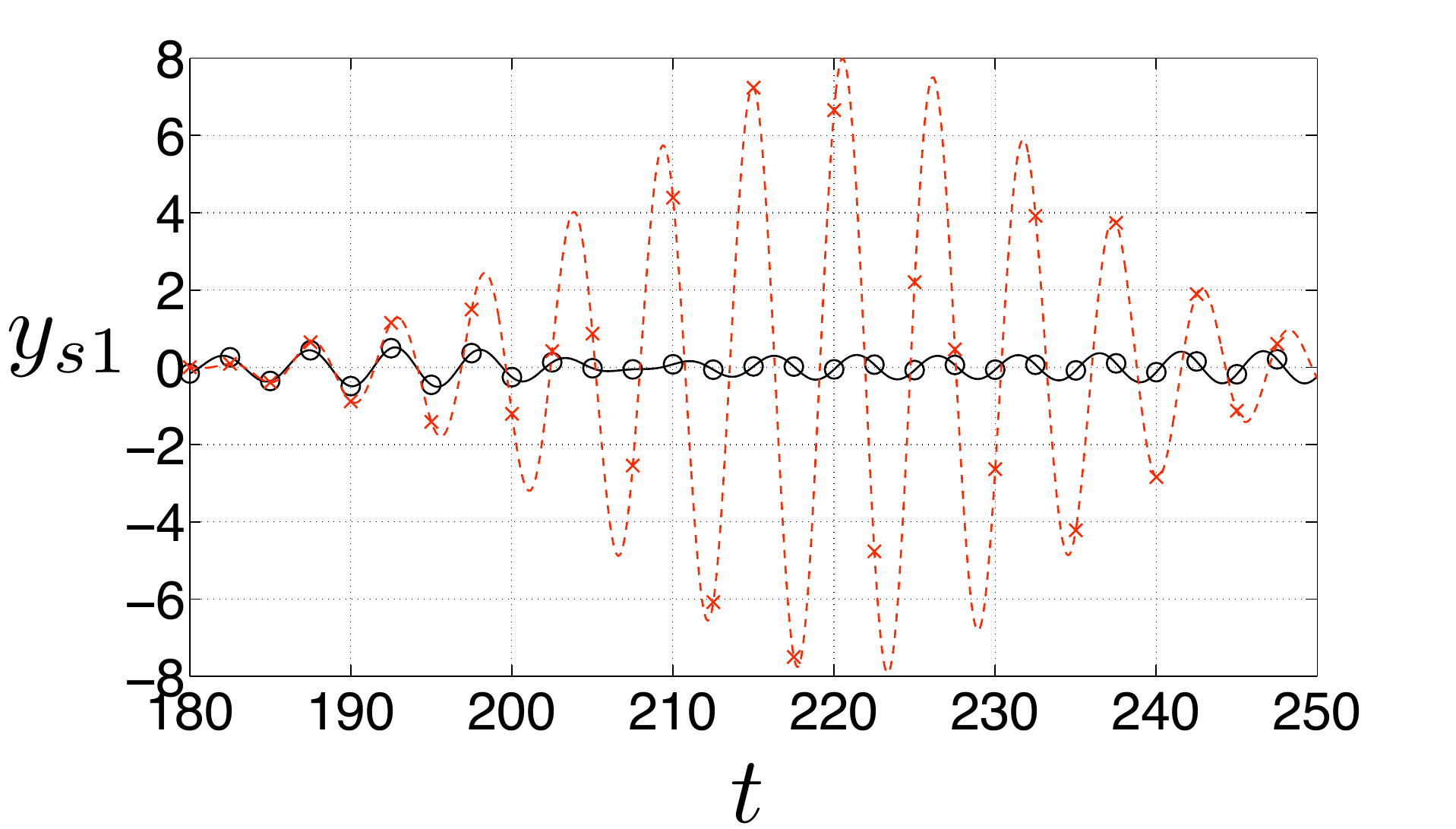}
\caption{Outputs of a system with full-state feedback control. The control gain is obtained using LQR, and the initial condition is that corresponding to~$t=180$ of the uncontrolled case plotted in Fig.~\ref{fig:cont_fullfback_lift}. Comparison of the outputs~$y_{u1}$ and~$y_{s1}$ of a 12-mode (2~unstable and 10~stable modes) reduced-order model ({\color{red}$\dashed,\times$}) with the projection of data from the full nonlinear simulation ($\solid,\circ$).} \label{fig:outputs_control_oproj4}
\end{figure}

We now use the same controller in the full nonlinear simulations and test the performance of the model for various perturbations of the steady state. A plot of the lift coefficient~$C_L$ vs.~time~$t$, with the control turned on at different times of the base simulation, is shown in Fig.~\ref{fig:cont_fullfback_lift}. The initial condition for the base case (no control) is the unstable steady state; eventually, small numerical errors excite the unstable modes and the flow transitions to periodic vortex shedding. In separate simulations, control is turned on at times~$t= 170, 180, 210$ corresponding to the base case. As the figure shows, the control is effective and is able to stabilize the steady state in each case, even when the flow exhibits strong vortex shedding. We remark that the latter two of these perturbations are large enough to be outside the range of validity of the linearized system, but the control is still effective, implying a large basin of attraction of the stabilized steady state. We also compare the output of the reduced-order model with the outputs of the nonlinear simulation; the plots are shown in Fig.~\ref{fig:outputs_control_oproj4}. The models perform well for the initial transients, but for longer times fail to capture the actual dynamics. This is not surprising as these perturbations are outside the range of validity of the linear models. For control purposes, it appears to be sufficient to capture the initial transients (approximately one period), during which the instability is suppressed to a great extent. We remark that one could possibly compute nonlinear models by projecting the full nonlinear equations onto the balancing modes, or enhance the model subspace by adding POD modes of vortex shedding and the shift modes as proposed by~\cite{Noack-jfm05} to account for the nonlinear terms.

Finally, we note that the reduced-order model~(\ref{rom_compact}) decouples the dynamics on the stable and unstable subspaces, and also, the dynamics on the unstable subspace can be computed only using the unstable eigen-bases~$\Phi_u$ and~$\Psi_u$. Thus, we could derive a control gain~$K \in \mathbb{R}^{1 \times n_u}$, based only on the two-dimensional unstable part of the model, such that the eigenvalues of~$(\widetilde{A}_u - \widetilde{B}_u K)$ are in the left half complex plane. That is, we can obtain a stabilizing controller {\em without} modeling the stable subspace dynamics. We have performed simulations to test such a controller and found that it also is capable of suppressing the periodic vortex shedding and thus results in a large basin of attraction for the stabilized steady state. The choice of weight matrices~$Q$ and~$R$ in the LQR cost~(\ref{lqr_cost}) needs to be different to obtain a comparable performance. However, as shown in the next section, it is essential to model the stable subspace dynamics to design a practical controller based on an observer that reconstructs the entire flow field using a few sensor measurements.

\subsection{Observer-based feedback control} \label{sec:observer}

The full-state feedback control of section~\ref{sec:fullstatefb} is not practical since it is not possible to measure the entire flow field. Here, we consider a more practical approach of measuring certain flow quantities at a small number of sensor locations. We assume that we can measure the velocities at the sensors shown in Fig.~\ref{fig:actuator_sensor}, in the near-wake of the plate. We remark that, even though these sensors are not experimentally realizable, they serve as a good testing ground for our models.

\begin{figure}
\centering
\includegraphics[height=1.5in]{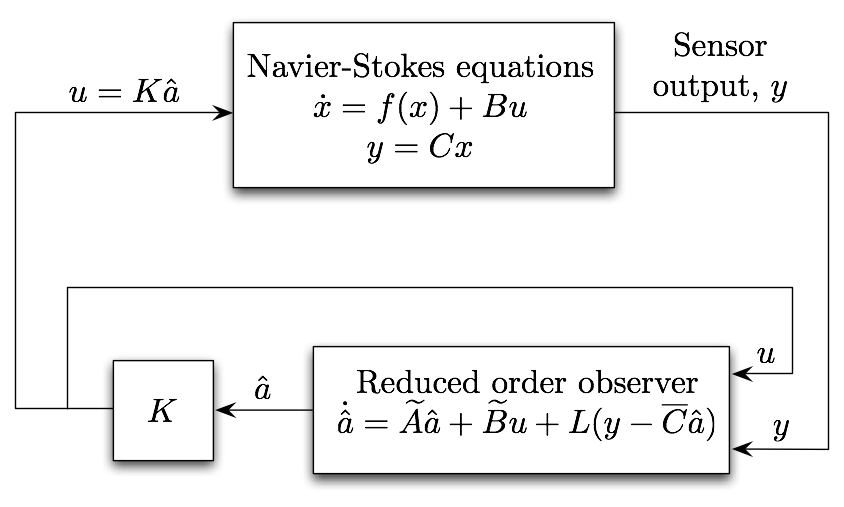}
\caption{Schematic of the implementation of observer-based feedback control in the nonlinear simulations. The control input~$u$ and the sensor measurements~$y$ are used as inputs to the observer, which reconstructs the reduced-order state~$\hat{a}$. This state is then multiplied by the gain~$K$, to obtain the control input~$u$. Both, the controller and observer gains~$K$ and~$L$ are computed based on the reduced-order model using LQR and LQG respectively.}
\label{fig:cartoon_obsv_control}
\end{figure}

\subsubsection{Observer design} \label{sec:obsv_design}
Using the reduced-order models derived as outlined in section~\ref{sec:models}, we design observers that dynamically estimate the entire flow field. Since the models~(\ref{rom_compact}) have a different output (projection onto certain modes), we modify the output equation so that it appropriately represents the sensor measurements. We replace the output equation of~(\ref{rom_compact}) with
\begin{align}
y &= M \begin{pmatrix} \widetilde{C}_u &  \widetilde{C}_s \end{pmatrix}  a  \overset{\text{def}}{=} \overline{C} a,
\label{observer_output}
\end{align}
where~$M \in \mathbb{R}^{s \times n}$ and~$s$ is the number of sensor measurements. The matrix~$M$~is sparse and extracts the values of the output of~(\ref{rom_compact}) at the sensor locations; thus, each row of~$M$ is filled with 0s except for the entry corresponding to a sensor measurement, which is~1. With the output equation~(\ref{observer_output}), we design an observer using a Linear Quadratic Gaussian (LQG) estimation. This method assumes that the errors in representing the state~$a$ and and the measurement~$y$ (due to the inaccuracies of the model) are stochastic Gaussian processes, and results in an estimate~$\hat{a}$ of the state~$a$ that is optimal in the sense that it minimizes the mean of the squared error; refer to~\cite{SkogPost-05} for details. We now discuss briefly our procedure for modeling these noises; consider the reduced-order model~(\ref{rom_compact}), but with process noise~$w$ and sensor noise~$v$ which enter the dynamics as follows:
\begin{align}
\dot{a} & = \widetilde{A} a + \widetilde{B} u + w \label{process_noise}\\
y & = \overline{C} a + v. \label{sensor_noise}
\end{align}
A key source of the process (state) noise~$w$ arises from model truncation, and second, from ignoring the nonlinear terms in the reduced-order model. The nonlinearity of the dynamics is important, for instance, when the model is used to suppress vortex shedding. A source of the sensor noise arises from two sources; first, the state~$x$ is approximated as a sum of a finite number of modes~(\ref{modal_exp}), and second, in the output projection step, the output is considered as a projection of the (approximated) state~$x$ onto a finite number of POD modes~(\ref{yrom_unstable_outproj}). Here, we approximate these two noises as Gaussian processes whose variances are
\begin{align}
Q &= E(w w^\ast), &w &= f(a_{meas})-\widetilde{A}a_{meas}, \label{process_noise_model} \\
R &= E(v v^\ast), &v &= y-\overline{C}a_{meas}, \label{sensor_noise_model}
\end{align}
and~$E(\cdot)$ gives the expected value. Here,~$f(\cdot)$ is the operator obtained by projecting the nonlinear Navier-Stokes equations~(\ref{navier1}) onto the balancing modes~$\Phi$, using the adjoint modes~$\Psi$. The state~$a_{meas}$ is obtained by projecting the snapshots, obtained from a representative simulation of the full nonlinear system, onto the balancing modes. The representative simulation we used here is the base case, with no control, shown in Fig.~\ref{fig:cont_fullfback_lift}, which includes the transient evolution from the steady state to periodic vortex shedding. The resulting estimator is of the form:
\begin{align}
\dot{\hat{a}} & = \widetilde{A} \hat{a} + \widetilde{B} u + L (y - \overline{C} \hat{a}), \label{lqg_observer}\\
\hat{y} & =  \overline{C} \hat{a},
\end{align}
where~$\hat{a}$ is the estimate of state~$a$, $\hat{y}$ is the estimated output, and~$L$~is the observer gain. The estimator is then used along with the full-state feedback controller designed in section~\ref{sec:fullstatefb} to determine the control input; a schematic is shown in Fig.~\ref{fig:cartoon_obsv_control}.

\begin{figure}
\centering
\includegraphics[height=2in]{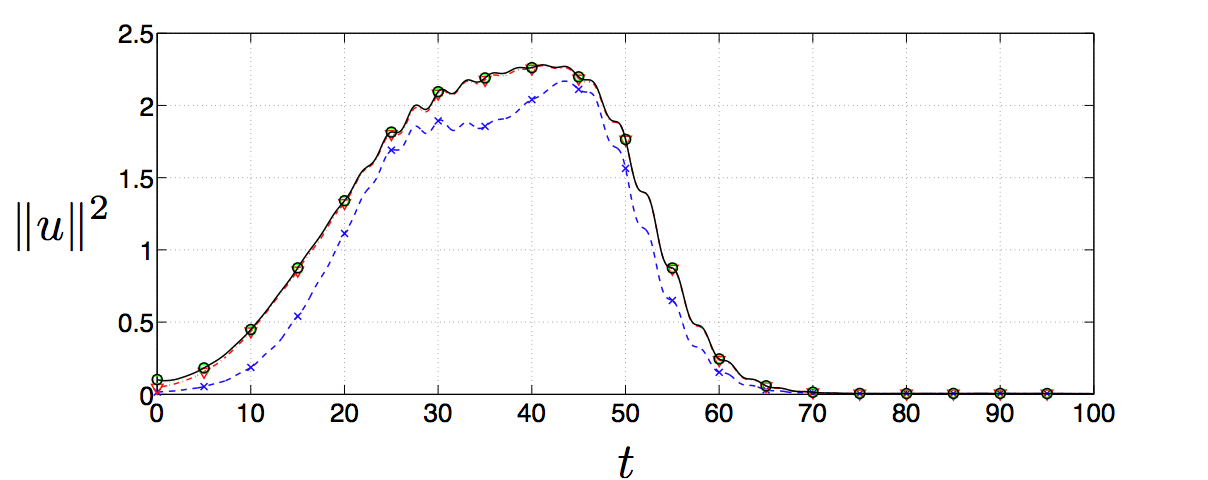}
\caption{The energy of the flow-field ($L^2-$norm) obtained from an impulse response~($\solid$,~$\circ$) of~(\ref{ss_stable_numerical}), and the energy captured by 4~({\color{blue}$\dashed$, $\times$}), 10~({\color{red}$\chndot$, $\triangledown$}), and 20~({\color{green}$\dotted$,~$\square$}) leading POD modes.}
\label{fig:l2norm}
\end{figure}

\begin{figure}
\centering
\includegraphics[height=2.7in]{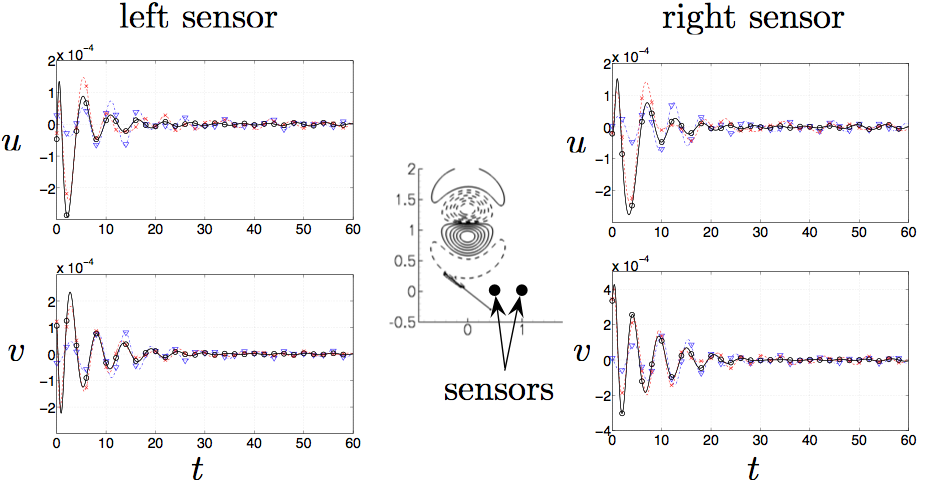}
\caption{Velocities at the sensor locations ($\solid$,~$\circ$), of an impulse response of~(\ref{ss_stable_numerical}), compared with the reconstruction using 10~({\color{blue}$\chndot$,~$\triangledown$}) and 20~({\color{red}$\dashed$,~$\times$}) leading POD modes.}
\label{fig:sensor_reconst}
\end{figure}

Since the observability Gramian corresponding to the pair~$(\widetilde{A}, \overline{C})$ is different from that for the pair~$(\widetilde{A}, \widetilde{C})$, the model~(\ref{rom_compact}) with the output represented by~(\ref{observer_output}), is not balanced. In principle, it is possible to construct reduced-order models of a system with the sensor measurements as the outputs, using the procedure of section~\ref{sec:approx_baltrunc_unstable}. Since the number of outputs will be typically small, the output projection step would not be required for such models. However, if this were done, the cost function~(\ref{lqr_cost}), based on total kinetic energy in the perturbation, would not be captured well by the model and the model would not be as effective for control design.

\begin{figure}
\centering
%{\includegraphics[width=4.5in]{compens_oproj20_lift}}
{\includegraphics[width=4.5in]{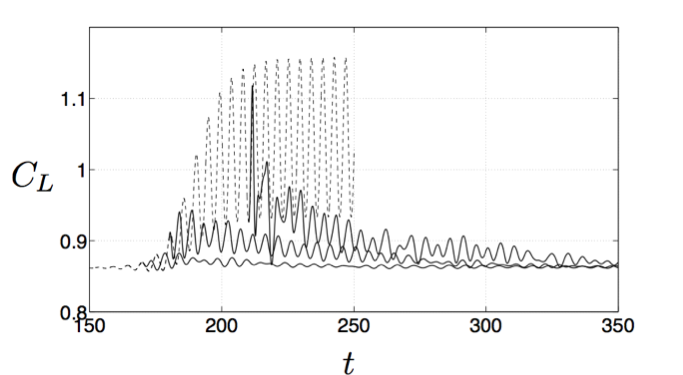}}
\caption{Lift-coefficient~$C_L$ vs. time~$t$, for estimator-based feedback control, with control turned on at different times in the base uncontrolled simulation. The base case~($\dashed$) is the same as in Fig.~\ref{fig:cont_fullfback_lift}, and the control is tested for different initial conditions, corresponding to~$t = 170, 180, 210$ of the base case~($\solid$). In both the cases, the controller stabilizes the flow to a small neighborhood of the steady state.} \label{fig:cont_obsvfback_lift}
\end{figure}

\begin{figure}
\centering
\includegraphics[height=1.4in]{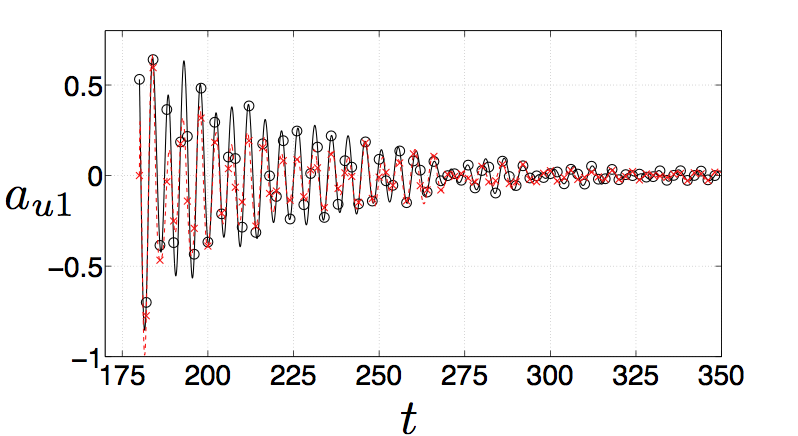}
\includegraphics[height=1.4in]{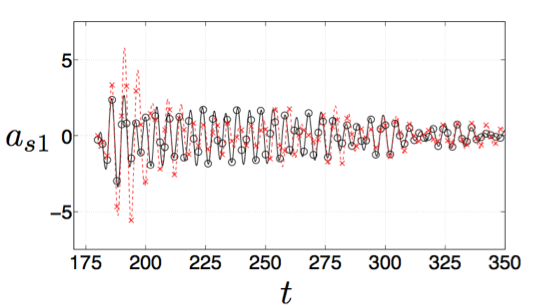} \\
\includegraphics[height=1.4in]{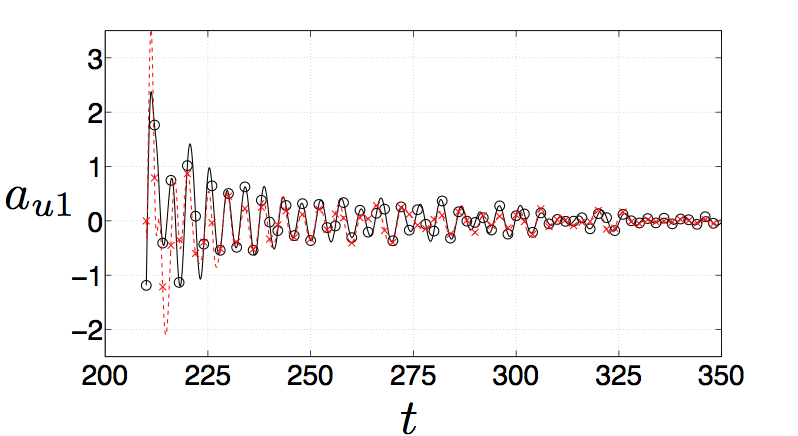}
\includegraphics[height=1.4in]{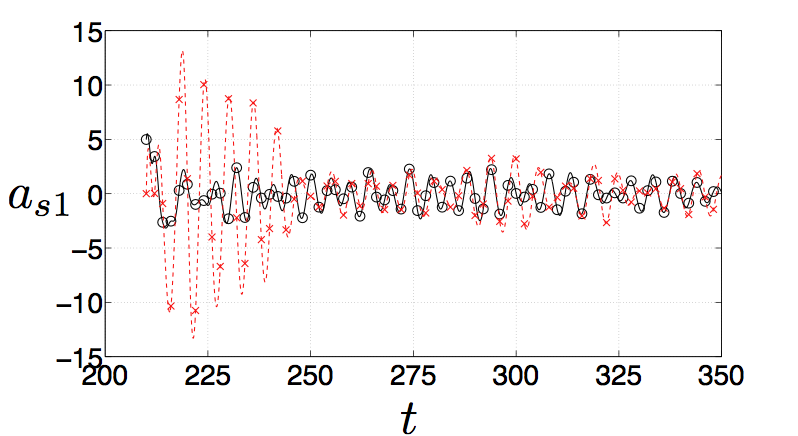}
\caption{States of the system with observer-based control; the states reconstructed~({\color{red}$\dashed$,~$\times$}) by a 22-mode observer quickly converge to the actual states~($\solid$,~$\circ$). The initial conditions used are those corresponding to~$t=180, 210$~(top and bottom) of the uncontrolled case shown in Fig.~\ref{fig:cont_obsvfback_lift}.}
\label{fig:states_obsvcontrol_oproj20}
\end{figure}

The models used in section~\ref{sec:fullstatefb} for full-state feedback were those of a system whose stable-subspace output was the velocity field projected onto the leading 4~POD modes. These 4~POD modes capture only about~$85\%$ of the energy, but the resulting models were effective in suppressing vortex shedding. However, for observer design, this representation of the output is inadequate, as the energy content of the flow at the sensor locations is very small, while the POD modes capture the energetically dominant modes. Hence, a greater number of POD modes is required to accurately represent the velocity at the sensor locations. The temporal evolution of the energy content of the flow, obtained from an impulse response of the system restricted to evolve on the stable subspace, is plotted in Figure~\ref{fig:l2norm}. Also plotted is the energy content of the same flow, but projected onto the leading 4, 10~and 20~POD modes; thus, a 4-mode projection leads to noticeable errors, while both 10-~and~20-mode projections accurately represent the energy. The velocity field at the sensor locations, reconstructed by 10~and 20~POD modes, plotted in Figure~\ref{fig:sensor_reconst}, shows that a 10-mode projection does not accurately represent the velocities at the sensor locations. Since 20 POD modes are sufficient to represent these velocities, we derive models using a {\em 20-mode output projection}, and use the same for observer design.

\subsubsection{Observer-based control}

The models obtained using the modified output~(\ref{observer_output}) are used to design dynamic observers based on the vertical~($v$-) velocity measurements at the sensor locations. A 22-mode reduced-order model, with 2~and 20~modes describing the dynamics on the unstable and stable subspaces respectively, is used to design a Kalman filter for producing an optimal estimate of the velocity field based on Gaussian approximations of error terms~(\ref{process_noise_model}, \ref{sensor_noise_model}). This estimate is then used along with reduced-order model controller to determine the control input, as shown in Fig.~\ref{fig:cartoon_obsv_control}. The results of this observer-based controller (or compensator) are shown in Figs.~\ref{fig:cont_obsvfback_lift},~\ref{fig:states_obsvcontrol_oproj20}. The compensator again stabilizes the unstable operating point, and furthermore, the observer reconstructs the reduced-order model states accurately. Initially, the observer has no information about the states (the initial condition is zero), but it quickly converges to and follows the actual states. There is a key difference from the full-state feedback control, that the  compensator does not exactly stabilize the unstable operating point but converges to its small neighborhood. The reason is that the observer design is based on the velocities at the  sensor locations projected onto the leading 20~POD modes, rather than the exact velocities at these locations.  These small errors enter into the observer's dynamics in the same way that sensor noise enters, resulting in small errors in the state estimate.

\section{Summary and discussion} \label{sec:discussion}

We presented an algorithm for developing reduced-order models of the input-output dynamics of high-dimensional linear unstable systems, extending the approximate balanced truncation method developed by~\cite{Rowley-ijbc05} for stable systems. We assumed that the dimension of the unstable eigenspace is small and the corresponding global eigenmodes can be numerically computed. The modeling procedure treats the dynamics on the unstable subspace exactly and obtains a reduced-order model of the dynamics on the stable subspace.

In a proof-of-concept study, the procedure was applied to control the 2-D~low-Reynolds-number flow past a flat plate at a large angle of attack~$\alpha$, where the natural flow state is periodic vortex shedding. We first performed a continuation study at~$\Rey = 100$ and computed the branch of steady states with~$\alpha$ varying from~$0$ to~$90^\circ$; we show that the flow undergoes a Hopf bifurcation from steady state to periodic shedding at ~$\alpha \approx 27^\circ$. We developed reduced-order models of the linearized dynamics at~$\alpha = 35^\circ$ actuated by a localized body force close to the leading edge of the plate. The outputs were considered to be the entire flow field, projected onto the unstable eigenmodes and the leading POD modes of the impulse response simulation (restricted to the stable subspace). We developed stabilizing controllers based on the reduced-order models to stabilize the unstable steady state and showed that the models agreed well with the actual simulations. We also included the controllers in the full nonlinear simulations, and showed that they had a large-enough basin of attraction to even suppress the vortex shedding. For such large perturbations, however, the model agreement with the full simulation was good only for short times. A natural step towards improving these models would be to project the full nonlinear equations onto the balancing modes to obtain nonlinear models. Alternately, the balanced models, which accurately capture the transient dynamics, could be combined with the POD-based models using shift-modes of~\cite{Noack-jfm05} which accurately capture vortex shedding and some of the transient dynamics. An interesting future direction is development of algorithms to compute {\em nonlinear} balanced models, for instance based on the theoretical work of~\cite{Scherpen-93}.

Instead of computing nonlinear models, here we pursued a step towards more practical controllers by considering an observer-based control design, in which the outputs were modified to be just two near-wake velocity measurements. The nonlinear terms in the equations, which our models do not capture, were treated as process noise, and the error in modeling the outputs was treated as sensor noise.
%For simplicity, both the noises were assumed to have a Gaussian distribution, though such an approach might not work in chaotic systems as pointed out by~\cite{KimBew-07}.
We designed a 22-mode reduced-order observer which reconstructed the flow field accurately, and along with the controllers, suppressed vortex shedding and stabilized the flow in a small neighborhood of the unstable steady state. We remark that the actuator and sensors considered here are not practically realizable, but the methodology presented here can be extended to a more practical actuation such as blowing and suction through the plate and measurements using surface pressure sensors. Furthermore, the choice of sensor locations in this study was ad hoc, and an interesting problem is of finding the optimal sensor locations, for a given actuator. The controllers present here are designed to operate at a fixed set of parameters~(such as~$\Rey$ and~$\alpha$), and it would be interesting to test their performance at off-design parameter values; the study on the performance of models of the linearized channel flow at off-design~$\Rey$ by~\cite{IlakRowley-pof08} shows promise in that direction.

A motivation for the choice of our model problem was to develop tools towards manipulating wakes of micro-air vehicles. Recently,~\cite{TairaCol-jfm08} performed a numerical study of flow past low-aspect-ratio plates, and a future direction we intend to undertake is to perform a detailed continuation study of the same flow to explore the existence and stabilization of high-lift unstable steady states in this 3-D flow.

\section{Acknowledgements}

The authors would like to thank Tim Colonius and Kunihiko Taira for their tremendous help in adopting their immersed boundary solver. The authors would also like to thank Ioannis G. Kevrekidis, Sung Joon Moon and Liang Qiao for their help with the timestepper-based steady-state analysis. This work was funded by the U.~S. Air Force Office of Scientific Research grant FA9550-05-1-0369 and this support is gratefully acknowledged.

\appendix

\section{Balancing transformation for unstable systems}
\label{sec:app:baltrunc}

Without loss of generality, a transformation~$T$ (and its inverse) that decouples the stable and unstable dynamics of~(\ref{ss}) can be written as:
\begin{align}
T = \begin{pmatrix} T_u & T_s \end{pmatrix}, \qquad T^{-1} = \begin{pmatrix} S_u^\ast \\ S_s^\ast \end{pmatrix},
\label{transf}
\end{align}
where the columns of~$T_u$ and~$T_s$ span the unstable and stable {\em right} eigenspaces of~$A$, while the columns of~$S_u$ and~$S_s$ span the unstable and stable {\em left} eigenspaces of~$A$. Further, these matrices are scaled such that~$S_u^\ast T_u = I_{n_u}$ and~$S_s^\ast T_s = I_{n_s}$. The transformation~(\ref{transf}) decouples the dynamics of~(\ref{ss}) as given in~(\ref{decouple_ss}) with the various matrices defined as follows:
\begin{align}
A_u &= S_u^\ast A T_u, \quad B_u = S_u^\ast B, \quad C_u = C T_u, \nonumber \\
A_s &= S_s^\ast A T_s, \quad B_s = S_s^\ast B, \quad C_s = C T_s.
\label{stable_def}
\end{align}
Using~(\ref{transf}) in~(\ref{gramians_general}), the Gramians of the original system~(\ref{ss}) are
\begin{align}
W_c & = T_u W_c^u T_u^\ast + T_s W_c^s T_s^\ast,  \nonumber \\
W_o & = S_u^\ast W_o^u S_u  + S_s^\ast W_o^s S_s,
\label{gramians_general_expand}
\end{align}
where, $W_c^s$ and~$W_o^s$ are the Gramians corresponding to the system defined by~$(A_s, B_s, C_s)$, while $W_c^u$ and~$W_o^u$ are the Gramians corresponding to the system defined by~$(-A_u, B_u, C_u)$. Let~$\widetilde{\Phi}_u \in \mathbb{R}^{n_u \times n_u}$ be the transformation that balances the Gramians~$W_c^u$ and~$W_o^u$, while~$\widetilde{\Phi}_s \in \mathbb{R}^{n_s \times n_s}$ be the transformation that balances~$W_c^s$ and~$W_o^s$. Then, it can be verified that the transformation that balances the Gramians~$W_c$ and~$W_o$ is given by
\begin{align}
\Phi = \begin{pmatrix} T_u \widetilde{\Phi}_u & T_s \widetilde{\Phi}_s \end{pmatrix} \overset{\text{def}}{=} \begin{pmatrix} \Phi_u & \Phi_s \end{pmatrix}.
\end{align}
Thus, the balancing transformation consists of two parts~$\Phi_u$ and~$\Phi_s$ which respectively balance the dynamics on the unstable and stable subspaces of~$A$. As per the technique of~\cite{ZhoSalWu-99}, a reduced-order model can be obtained by truncating the columns of~$\Phi$ that correspond to the relatively uncontrollable and unobservable states. As we will show now, the algorithm outlined in section~\ref{sec:approx_baltrunc_unstable} essentially computes the leading columns of~$\Phi_s$ (and the corresponding rows of its inverse).
%Similarly, it can be shown that the inverse transformation is given by~$\Psi = \begin{pmatrix} S_u \widetilde{\Psi}_u & S_s \widetilde{\Psi}_s \end{pmatrix}.$
We show that the controllability Gramian of the stable dynamics of~(\ref{ss}), which are defined by~(\ref{ss_stable}), is the same as the ``stable'' part of the Gramian defined in~(\ref{gramians_general_expand}). Note that using~(\ref{transf}) and the definition~(\ref{proj_stable}), the projection operator~$\mathbb{P}_s$ can be written as
\begin{align}
\mathbb{P}_s = I - T_u S_u^\ast = T_s S_s^\ast.
\label{proj_stable_alt}
\end{align}
Using the definition~(\ref{gramians_stable}), the controllability Gramian of~(\ref{ss_stable}) is
\begin{align}
\widetilde{W}_c^s &= \int_0^\infty \, e^{\mathbb{P}_s A t}  (\mathbb{P}_s B) \, (\mathbb{P}_s B)^\ast e^{ (\mathbb{P}_s A)^\ast t} \, \, dt \nonumber \\
& =  \int_0^\infty \, T_s \, \,  e^{ S_s^\ast A T_s t}  \, \, S_s^\ast B \,\, B^\ast S_s \, \,  e^{ T_s^\ast A^\ast S_s t} \, \, T_s^\ast \, \, dt   \qquad \mbox{using equation~(\ref{proj_stable_alt})} \nonumber \\
& =    T_s \Big(\, \, \int_0^\infty \, e^{A_s t} \, \, B_s \, B_s^\ast \, \, e^{A_s^\ast t} \, \, dt  \, \, \Big) T_s^\ast \qquad \qquad \mbox{using equation~(\ref{stable_def})} \nonumber \\
& = T_s W_c^s T_s^\ast,
\end{align}
which is the same as the stable part of~$W_c$. Similarly, it can be shown that the observability Gramian~$\widetilde{W}_o^s$ of~(\ref{ss_stable}) is the same as the ``stable'' part of the observability Gramian~$W_o$:
\begin{align}
\widetilde{W}_o^s &= \int_0^\infty \, e^{\mathbb{P}_s^\ast A^\ast t} ( \mathbb{P}_s^\ast C^\ast) \, (\mathbb{P}_s^\ast C^\ast)^\ast e^{ (\mathbb{P}_s^\ast A^\ast)^\ast t} \, \, dt = S_s^\ast W_o^s S_s.
\end{align}
Thus, balancing the Gramians~$\widetilde{W}_o^s$ and~$\widetilde{W}_c^s$ is identical to balancing the parts of the Gramians~$W_c$ and~$W_o$ of the original system~(\ref{ss}) that are related to the dynamics on the stable subspace of~$A$.

\section{Derivation of the adjoint equations}
\label{sec:app:adjoint}

In this appendix, we derive the adjoint of the linearized semi-discrete equations~(\ref{linear},\ref{lin_constraint}).
Let~$(\zeta, \psi)$ be the weighting functions corresponding to~$(\gamma, \tilde{f})$. Then, using the inner product defined in equation~(\ref{ip}), the weak form of~(\ref{linear},\ref{lin_constraint}) is:
\begin{align}
\int_0^T \int_\Omega \, \zeta \cdot (C^T C)^{-1} & \,\Big(\frac{d\gamma}{dt} + C^T E^T \tilde{f} + \beta C^T C \gamma - C^T {\cal N}_L(\gamma_0) \gamma \Big) \, dx \, dt \nonumber \\
& + \int_0^T \int_\Omega \, \psi \cdot ECs \, dx \,dt = 0.
\label{weak1}
\end{align}
Integrating by parts with respect to~$t$ and rearranging terms,
\begin{align}
\int_0^T \int_\Omega \,\, \gamma \cdot  \, & \Big( - (C^T C)^{-1}\frac{d\zeta}{dt} + (C^T C)^{-1} C^T E^T \psi + \beta \zeta - \big((C^T C)^{-1} C^T {\cal N}_L(\gamma_0) \big)^T \zeta \Big) \, dx \, dt \nonumber \\
& + \int_0^T \int_\Omega \,\,\tilde{f} \cdot \Big(EC (C^T C)^{-1} \zeta\Big) \, dx \, dt + \ip<\gamma,\zeta> \Big|_0^T = 0.
\label{weak2}
\end{align}
For linearization about stable steady states,~$\gamma \rightarrow 0$, as~$T \rightarrow \infty$, and if the adjoint equations are integrated backwards in time, $\zeta(t=0) \rightarrow 0$. So, the last term on the left hand side of equation~(\ref{weak2}) vanishes identically. If equation~(\ref{weak2}) is to hold for all values of~$\gamma$ and~$\tilde{f}$, we get the following adjoint equations hold:
\begin{align}
-\frac{d\zeta}{dt} + C^T E^T \psi & = - \beta C^T C \zeta + (C^T C){\cal N}_L(\gamma_0)^T q_a , \label{adjoint1} \\
EC \xi & = 0, \label{adjoint2}
\end{align}
where~$\xi = (C^T C)^{-1} \zeta$ and $q_a = C \xi$ can be thought of as the weighting functions corresponding to the streamfunction~$s$ and the flux~$q$ respectively. Now, equations~(\ref{adjoint1},\ref{adjoint2}) have the same form as~(\ref{linear},\ref{lin_constraint}) except for the nonlinear term. Thus, the same time-integration scheme can be used for both, with the appropriate (linearized) nonlinear terms.

\bibliographystyle{jfm}
%\bibliography{jabbrv,master}

\end{document}